\documentclass[12pt]{article}
\usepackage[amsmath]{my_e-jc}

\usepackage{amsmath}
\usepackage{amsthm}
               
\usepackage{graphicx}     
\usepackage{hyperref} 
\usepackage{url}
\usepackage{amsfonts} 

\usepackage{multicol}
\usepackage{textcomp}
\usepackage{ucs} 
\usepackage{comment}
\usepackage{rotating}

\usepackage{graphicx}
\usepackage{yhmath}
\usepackage{mathdots}
\usepackage{MnSymbol}
\usepackage{adjustbox}
\usepackage[graphicx]{realboxes}
\usepackage{graphicx}
\usepackage{float,lscape}
\usepackage{bigstrut}
\usepackage{diagbox}
\usepackage{colortbl}
\usepackage{tabularx}
\usepackage{verbatim}
\usepackage{algorithm}
\usepackage{algpseudocode}
\usepackage{tabularray}

\usepackage{gensymb}

\usepackage{pifont}

\allowdisplaybreaks

\makeatletter
\@addtoreset{footnote}{page}
\makeatother

\dateline{Sep 21, 2025}{TBD}{TBD}

\MSC{05A19}

\Copyright{The authors. Released under the CC BY-ND license (International 4.0).}

\title{Interpreting a Sum of Third Powers by using\\ 
       a Geometric Assembly in Four Dimensions}
       
\author{Joseph Alfano\authornote{1}   \and
        Emily Armstrong\authornote{2} \and
        Vincent DeNolf\authornote{3}  \and
        Jonah Sagarin\authornote{4} }

\authortext{1}{Department of Mathematics and Computer Science, 
               Assumption University, Worcester, Massachusetts, USA 
               (\email{jalfano@assumption.edu}).}  
\authortext{2}{Department of Mathematics, 
               University of Pittsburgh, Pittsburgh, Pennsylvania, USA 
               (\email{ema203@pitt.edu}).}               
\authortext{3}{University of Colorado, Boulder, Colorado, USA 
               (\email{vincent.denolf@colorado.edu}).}
\authortext{4}{Union College, Schenectady, New York, USA 
               (\email{sagarinj@union.edu}).}

\begin{document}

\maketitle

\begin{abstract} 
 We present new combinatorial proofs of Nicomachus's Theorem 
 for the sum of the third powers of the first n natural numbers.
 The key step is that we define a 4-dimensional block which comprises 
 unit hyper-cubes.
 In our first proof we assemble 4 copies of this block
 to construct a rectangular solid.
 In our second proof we partition the block into two parts, 
 then map and reassemble them into a solid 
 whose shape is a step triangle along two coordinate axes, 
 and is likewise on the other two.
 For corollaries we present a q-analogue of this identity using taxicab distance,
 and we interpret q-analogues from 4 different (sets of) authors,
 using taxicab distances from different starting points.
\end{abstract}

\section{Introduction}

Identities for evaluating the sum 
$1^p + 2^p + \ldots + n^p$
for a positive integer $p,$ such as 
\begin{align}  
  1^2 + 2^2 + \cdots + n^2
  &= 
  \tfrac{1}{3} n (n+1) \left(n+\tfrac{1}{2}\right)\\
  1^3 + 2^3 + \cdots + n^3
  &= \tfrac{1}{4} n^2(n+1)^2,\hspace{3.50em}\text{Nicomachus's Theorem}
\end{align} 
are a popular inspiration for combinatorial interpretation.

M. K. Siu proved identity (1) by defining a 3-dimensional block which comprises
{${1^2 + 2^2 + \cdots + n^2}$} unit cubes,
then assembling 3 copies of this block
and then relocating a half-unit high subcomponent, 
to construct a rectangular solid with side lengths 
$n$, $n+1$, and $n + \tfrac{1}{2}$. 
More precisely, Siu's construction interprets an identity
which is algebraically equivalent to (1): 
$\hspace{0.1250em}(1^2 + 2^2 + \ldots + n^2)\cdot 3 = 
               n(n+1)(n+\tfrac{1}{2}).\hspace{0.1250em}$
(See \cite{example3_siu}.)

And several authors have proved identity (2) using geometric counting arguments:
by assembling unit squares 
  (see \cite{example17_berendonk,
             example5_cupillari,
             example9_flores,
             example1_golomb,
             example2_love}) 
and/or triangles in 2 dimensions 
  (see \cite{example15_laosinchai,
             example7_schrage}),
or unit cubes in 3 dimensions 
  (see \cite{example4_fry,
             example13_kanim,
             example6_nelsen, 
             example18_stevanovic_etal}).

In this article, we present a proof of identity (2)
by extending Siu's construction to 4 dimensions. 
We define a 4-dimensional block which comprises 
$1^3 + 2^3 + \cdots + n^3$ unit hyper-cubes
(which we call \emph{unit cubes}), 
then we assemble 4 copies of this block and 
then relocate a unit-high subcomponent,
to construct a rectangular solid with side lengths $n$, $n$, $n+1$, $n+1$.
More precisely, this construction interprets an identity which is
algebraically equivalent to (2): 
$\hspace{0.1250em}(1^3 + 2^3 + \ldots + n^3) \cdot 4 = n^2 (n+1)^2.$

And we present a second proof of (2),
in which we take a single copy of the 4-dimensional block 
and partition it into two parts, 
then map and reassemble them into a solid 
whose shape along two coordinate axes is a step triangle
comprising $n(n+1)/2$ unit squares, 
and along the other two axes is also a step triangle 
comprising $n(n+1)/2$ unit squares.
This interprets (2) more directly, 
without multiplying each side by 4. 
\\

Identity (2) is also an inspiration for constructing $q$-analogues. 
Several authors have found solutions to this question, 
by using for interpretation:
integer partitions (Garrett and Hummel~\cite{example10_garrett_etal}); 
telescoping argument (Warnaar~\cite{example11_warnaar});               
lattice of subspaces (Zhao and Feng~\cite{example14_zhao_etal});      
and overpartitions (Forster, Garrett, Jacobsen, 
                    and Wood~\cite{example16_forster_etal}).           

At the end of Section 2 of this article
which presents our four-block construction, 
we introduce a $q$-analogue of identity (2), 
in its algebraically equivalent form.
We use the taxicab distance from a certain start point 
to each unit cube in our final fitted assembly:
one side of our equation records the taxicab distance 
to the unit cubes which constitute 
the entire rectangular solid;
and the other side records the taxicab distance 
to the unit cubes in each of the 4 copies of our block.

In Section 3 which presents our one-block construction, 
we also show that 
for each of the four $q$-analogues from the above authors, 
there is an interpretation which applies taxicab distance 
to our single-block proof of (2), each from a different start point.
One side of their equation records the 
taxicab distance to the unit cubes which constitute 
the reassembled solid;
and the other side records the taxicab distance to the unit cubes in 
each of the 2 mapped parts from the partition of our 
single 4-dimensional block.

\section{The four-block construction}
 
In this section we interpret identity (2), 
in the algebraically equivalent form
\begin{equation*}
(1^3 + 2^3 + \cdots + n^3) \cdot 4 = n^2 (n+1)^2,
\end{equation*}
using a geometric assembly of four 4-dimensional blocks.
    
We work in $\mathbb{R}^4$: label the coordinate axes $x,y,z,w.$
We fix these axes to be pairwise orthogonal,
but do not require that their orientation obey a right-hand rule.
We will find it helpful for visualization in the figures
to depict the $z$-axis oriented in the opposite direction
of the right-hand rule with respect to the $x$- and $y$-axes.
\\

Our fundamental object is the unit hyper-cube
\begin{equation}
\{(x,y,z,w) \mid 0 \le x,y,z,w \le 1\},
\end{equation}
and more generally,
its image under translation by any $(a-1,b-1,c-1,d-1) \in \mathbb{R}^4$:
\begin{equation}
\{(x,y,z,w) \mid a-1 \le x \le a;\thickspace 
                   b-1 \le y \le b;\thickspace
                   c-1 \le z \le c;\thickspace
                   d-1 \le w \le d\}.
\end{equation}
For brevity we call this object in (4) a \emph{unit cube,}
and omit the prefix ``hyper-"; there should be no confusion.
We will use the name $c_u$ to refer to
any arbitrarily selected such unit cube.
\\

In our assembly algorithm, we will apply to each unit cube $c_u$ 
one or more isometry mappings $f,$ i.e. some composition of 
translations, rotations, and/or reflections in $\mathbb{R}^4$.

For this article  
we restrict our selection of rotations and reflections $f ,$ 
to those which satisfy 
that the positions of our $c_u$ and $f(c_u)$ will always have the form (4). 
That is, for each of our four coordinate axes $x_j\thinspace ,$
the image of that axis under $f$ will be parallel 
to one of the four coordinate axes of the pre-image. 

Also for unit cube $c_u\thinspace ,$ in its form (4)
we restrict the values of $a,b,c,d$ to be integers.
And we restrict our selection of isometry mappings $f$ 
to those which satisfy
that the image $f(c_u)$ has integer valued endpoints 
on its unit interval for each coordinate axis. 

So as a convention to simplify notation, 
we assign to each unit cube $c_u$ that satisfies (4) 
a \emph{location label} by writing $(x,y,z,w) = (a,b,c,d).$
That is, for the unit cube location label
we write only the right endpoint of the unit interval 
which $c_u$ occupies, for each coordinate axis.

And when we write an individual location label coordinate such as $w=d,$ 
we intend it to signify the unit-length interval $d-1 \le w \le d$
(unless explicitly stated otherwise).

This convention will simplify the visualization in the figures:
for a given unit cube $c_u$ whose position is as stated in (4),
we depict that 4-dimensional cube 
by sketching the 3-dimensional cube
$\{(x, y, z) \mid a-1 \le x \le a;\thickspace 
                  b-1 \le y \le b;\thickspace 
                  c-1 \le z \le c\}$
at the coordinate location labeled $w=d.$

\subsection{Defining the blocks}

Fix an integer $n \ge 2.$
Here we construct four blocks 
(as in children's toys), which we call $A,B,C,D.$
Each block comprises $1^3 + 2^3 + \ldots + n^3$ unit cubes.

\begin{definition}
We let the \emph{blocks} $A,B,C,D$ each be 
the union of the unit cubes 
whose positions have the form (4), 
whose location labels $(x, y, z, w)$ satisfy the restrictions: 
\begin{alignat*}{2}
A & = \bigcup_{i=1}^n \,
      \bigcup_{\substack{
               z=i\\
               x, y, w \in \{1, 2, \ldots, i\}
                }
                }           
    &&( \textrm{unit cube at } (x,y,z,w) )
    \\
B & = \bigcup_{i=1}^n \,
      \bigcup_{\substack{
               y=i+1\\
               x, z, w \in \{1, 2, \ldots, i\}
                }
                }           
    &&( \textrm{unit cube at } (x,y,z,w) )
    \\    
C & = \bigcup_{i=1}^n \,
      \bigcup_{\substack{
               x=i+1\\
               y \in \{2, 3, \ldots, i+1\}\\
               z, w \in \{1, 2, \ldots, i\}
                }
                }           
    &&( \textrm{unit cube at } (x,y,z,w) )
    \\    
\hspace*{7.3800em}D & = \bigcup_{i=1}^n \,
      \bigcup_{\substack{
               w=i\\
               x,y \in \{1, 2, \ldots, i\}\\
               z   \in \{0, 1, \ldots, i-1\}
                }
                }           
    &&( \textrm{unit cube at } (x,y,z,w) ).\hspace*{11.6800em}
    {\rule{2pt}{1em}}   
\end{alignat*}    
\end{definition}


\noindent
See Figure 1 (at its step 1), for illustration of these blocks. 
There we have depicted them in an exploded-view positioning
within each unit interval in $w$, for clearer distinction.  
And we added visualization aids: 
\begin{itemize}
\item 
For each block, we added dashed lines:
to illustrate where two unit cubes $c_u$ and $c_u^{\prime}$
of a given block which appear to be non-adjacent are actually adjacent.
We did this where $c_u$ and $c_u^{\prime}$
have adjacent location label coordinates on the $w$-axis, 
$w=d$ and $w=d+1$ for some $d,$
and their location label coordinates in $(x,y,z)$ are identical.

\item 
For each block, we added a colored ball at the center of certain unit cubes: 
blue color
to highlight the ``indicator axis" for that block [to be defined below]; and
other colors
to highlight each of the non-indicator axes for that block.

\item 
For each block, we added arrows: 
to illustrate how these blocks assemble together 
in the ``first fitting" step;
and the path of the subcomponent which will get mapped in the 
``second fitting" step.
[Each of these fittings will be defined below.]
\end{itemize}
\begin{figure}
    \centering
    \includegraphics[scale=0.81250] {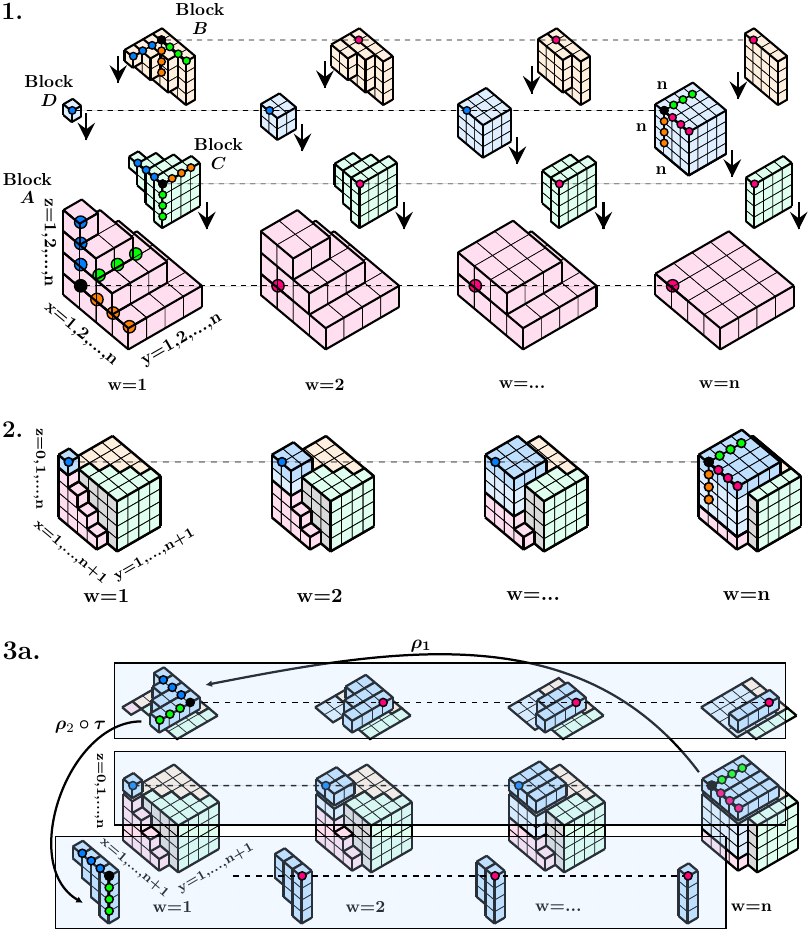}
    \caption{
    (1.) We interpret the sum $(1^3 + 2^3 + \ldots + n^3) \cdot 4$ using four 
    $4$-dimensional blocks, $A, B, C, D,$ here depicted for $n=4 .$ 
    Each block is a union of $\sum_{i=1}^{n} i^3$ unit hyper-cubes 
    (we call \emph{unit cubes}). 
    Within each block's defined location, its unit cubes each occupy 
    a unit interval with integer endpoints in the coordinates $x,y,z,w .$
    We label a unit cube's location by writing $(x,y,z,w) = (a,b,c,d) ,$ 
    recording the right endpoint of its unit interval on each axis. 
    E.g. our $w=1$ signifies the real interval $[0,1]$ for $w .$ 
    So each $4$-dimensional cube is depicted here 
    as a $3$-dimensional cube in $x,y,z ,$ at some $w=d .$
    \\
    \hspace*{1.0em} Dashed lines indicate that 
    wherever a block occupies two adjacent unit intervals in $w ,$ 
    then pairs of corresponding unit cubes are adjacent. And colored balls within 
    the unit cubes highlight that 
    our four blocks are equivalent, under rotations and translations.
    $D$'s equivalence with $C$ will follow from mappings to be depicted in Figure 2.
    \\
    (2.) The first fitting: our assembled blocks fill a 
    rectangular region, the union of the unit cube sized locations at 
    $\{ (x,y,z,w) \mid x,y\in \{1, 2, \ldots, n+1 \}; 
                       z,w\in \{1, 2, \ldots, n   \} \} ,$
    which we call $R_3(n) .$ 
    But there are overhang unit cubes at $z=0 ,$ and vacant locations at $y=1 .$
    \\
    (3a.) The second fitting: we map the overhang to fill the vacancy, by a sequence 
    of rotations and translations. We selected the 
    $\rho_1$ and $\rho_2 \circ \tau$ for this figure, to facilitate 
    visualization of the fitting through isometry mappings. Continued at Figure 2.
    }             
\end{figure}

\begin{figure}
    \centering
    \includegraphics[scale=0.8820] {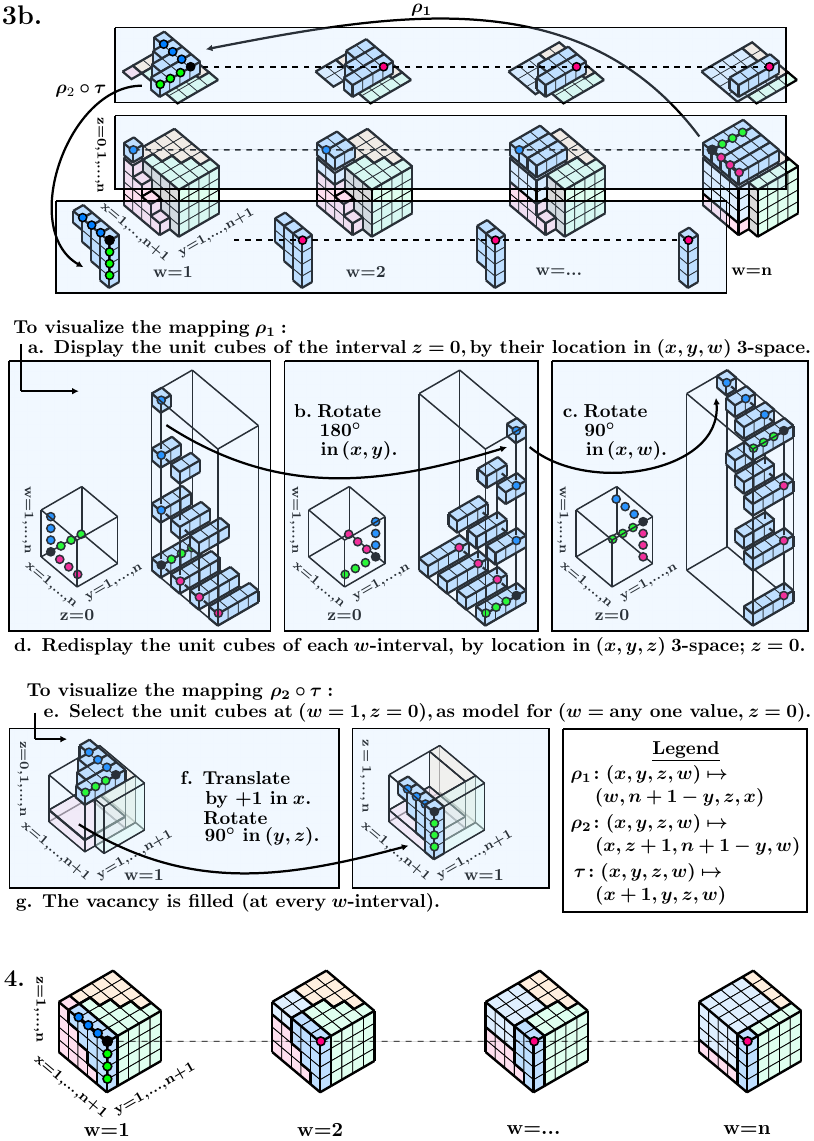}
        \caption{
    (3b.) Our isometry mappings $\rho_1, \rho_2, \tau$ are defined in the legend.
    To help visualize the rotation $\rho_1$:
    (a.) we change our perspective of the unit cubes at $z=0 ,$
    to display their location  in the $(x,y,w)$ $3$-space;
    (b, c.) rotate this set of unit cubes in the $(x,y)$ plane, and in the
    $(x,w)$ plane;
    (d.) change perspective back, to display these unit cubes at each unit 
    $w$-interval by their location in $(x,y,z)$ $3$-space, at $z=0 .$
    To help visualize the translation and rotation $\rho_2 \circ \tau$:
    (e, f, g.) we display their actions at $w=1 ,$ as representative for
    $w = $ any one value (at $z=0$).
    \\
    (4.) The final assembly comprises all of the 
    $(1^3 + 2^3 + \ldots + n^3) \cdot 4$ initial unit hyper-cubes,
    and fills completely the rectangular region $R_3(n) .$
    Its volume equals $n^2 (n+1)^2 .$
    }
\end{figure}

For nomenclature: in 3-dimensional geometry, a ``square pyramid" 
is a polyhedron formed by connecting a square base and an apex point, 
located outside the plane that contains the square. 
Each base edge and the apex form a triangle. 
If the apex is not located directly above the center of the square, 
we have an ``oblique square pyramid". 

In 4-dimensional geometry, a ``cubic pyramid" 
is a polytope formed by connecting a cube base and an apex point, 
located outside the 3-dimensional space that contains the cube.
Each square face of the cube and the apex form a square pyramid.

Our blocks $A,B,C,D$ are defined in a discrete manner.
For example in block $A ,$
its 3-dimensional subcomponent at $w=1 $
has shape which resembles a ``step pyramid" architectural structure, 
not pyramid with smooth triangular faces.
We will call the shape of this 3-dimensional subcomponent of block $A$ 
an \emph{(oblique) square step pyramid}.

Now we will call the shape of each of our entire 4-dimensional blocks $A,B,C,D$
an \emph{(oblique) cubic step pyramid}.
And for a 2-dimensional face 
like that of block $A$'s subcomponent at $w=1$ which lies in the $(x,z)$-plane,
that comprises $1 + 2 + \ldots + n$ unit squares,
we will call the shape of this face 
an \emph{(oblique) step triangle}.
\\

In our assembly algorithm, we will apply to each block 
(or subcomponent of a block) one or more isometry mappings $f$. 
The action of $f$ on a block is defined by 
the action of $f$ on its constituent unit cubes.

When we speak of the location in $\mathbb{R}^4$
of a given block, e.g. before or after mapping by 
a transformation $f,$ we intend that
we are speaking of the locations of its unit cubes.

The \emph{defined location} of each of our blocks $A,B,C,D$
is its location as stated in its definition above.

\subsection{Block equivalence}

Now observe that for each of these blocks, 
there is a unique translation mapping that sends the block,
in its defined location, 
to a position where all of the location label coordinates
of all of its unit cubes lie in $\{1, 2, \ldots, n\}.$

\begin{definition}
For each block $A, B, C, D,$ 
the \emph{canonical location} is that location 
(obtained by translating the block from its defined location) 
where all the location label coordinates
of all of its unit cubes lie in $\{1, 2, \ldots, n\}.$
\hfill {\rule{2pt}{1em}}
\end{definition}

\noindent
One checks that: 

Block $A$ is mapped from its defined location to %
its canonical location by the translation
which leaves each coordinate value unchanged.

Block $B$ is mapped from its defined location to %
its canonical location by the translation 
$\tau_{y-1}:(x,y,z,w) \mapsto(x,y-1,z,w).$

Block $C$ is mapped from its defined location to 
its canonical location by the translation 
$\tau_{x-1, y-1}:(x,y,z,w) \mapsto(x-1,y-1,z,w).$

Block $D$ is mapped from its defined location to 
canonical location by the translation 
$\tau_{z+1}:(x,y,z,w) \mapsto(x,y,z+1,w).$
\\

We will utilize both the defined location and the canonical location
for each block, in our assembly algorithm. 
The first-fitting step will assemble the blocks 
each in their defined location.
And the second-fitting step will map a subcomponent of one block onto a 
subcomponent of another block. 
The fact that such a bijective isometry mapping exists, is due to 
properties of blocks (``equivalence" and ``replication") 
which will be seen most easily when each is in canonical location.
\\

\begin{definition}
For each block $A, B, C, D,$ the \emph{indicator axis} 
is the unique coordinate axis $x_j$ which satisfies that: 
when the block is in canonical location, 
then for every one of its unit cubes $c_u,$
the largest coordinate value in the location label of $c_u$ 
is attained on the $x_j$ axis. \hfill {\rule{2pt}{1em}} 
\end{definition}

\noindent 
The indicator axis for the block $A$ (respectively $B, C, D$)
is the axis $z$ (respectively $y, x, w$).
\\

We will utilize the indicator axis for each block,
when we prove (and depict in the figures)
that our four blocks are mutually equivalent
under isometry transformations.\\

\begin{theorem}[Block equivalence]\label{Thm:2.1}
The blocks $A, B, C, D,$ when each in its defined location,
are pairwise equivalent under 
compositions of rotations and translations in $\mathbb{R}^4$.
\end{theorem}

\begin{proof}\phantom{\qedhere}
\hspace*{-2.0em}
Block $B$ is the image of $A$
under a rotation and translation,
$\tau_{y+1} \circ \rho_{y, z, x}$:
\begin{alignat*}{6}
\hspace{-0.72em}&\rho_{y, z, x}\hspace{0.90em}&&: (x, y, z, w) 
\mapsto (y&&,z&&,x&&,w&&)\\ 
\hspace{-0.72em}&\tau_{y+1}\hspace{0.90em}&&: (x, y, z, w) 
\mapsto (x&&,y+1&&,z&&,w&&). 
\end{alignat*}
And block $C$ is the image of $A$
under a rotation and translation,
$\tau_{x+1, y+1} \circ \rho_{z, x, y}$:
\begin{alignat*}{6}
&\rho_{z, x, y} 
&&: (x, y, z, w) 
\mapsto (z&&,x&&,y&&,w&&)\\
&\tau_{x+1, y+1} 
&&: (x, y, z, w) 
\mapsto (x+1&&,y+1&&,z&&,w&&). 
\end{alignat*}
Surprisingly, block $D$ is the image of $A$
under a rotation and translation,
$\tau_{z-1} \circ \rho_{w, x, z}$:
\begin{alignat*}{6}
\hspace{-0.69em}&\rho_{w, x, z}\hspace{0.78em}&&: (x, y, z, w) 
\mapsto (w&&,y&&,x&&,z&&)\\
\hspace{-0.69em}&\tau_{z-1}\hspace{0.78em}&&: (x, y, z, w) 
\mapsto (x&&,y&&,z-1&&,w&&).\tag*{\qed}
\end{alignat*}
\end{proof}

In other words: the blocks $A,B,C,D,$ 
when each is in canonical location, 
are pairwise equivalent under the set of rotations in $\mathbb{R}^4.$
\\

See Figure 1 (at its step 1), for illustration of the block equivalence.
For all four blocks, we added colored balls with 
same single color (blue) along the indicator axis, 
to aid the visualization that these blocks 
are equivalent under rotations and translations.

For our proof of the block equivalence from $A$ to $B,$ 
we selected the particular rotation $\rho_{y,z,x}$
for its intuitive clarity in Figure 1: 
it keeps the $w$-coordinate fixed,
so this rotation acts within the familiar $(x,y,z)$ 
3-dimensional space,
for each unit interval $w=d.$
But in this proof it would have sufficed for us 
to select \emph{any} permutation of the coordinate axes
which satisfies that: 
it maps the indicator axis $z$ of block $A$
to the indicator axis $y$ of block $B;$ 
and the permutation has even parity so it is a rotation. 
And likewise for our selection of the 
rotation $\rho_{z,x,y}$ from $A$ to $C.$

The equivalence of the block $D$ is the least intuitive visually, 
at this moment in the presentation.
We will explain below, our rationale 
for selecting the particular rotation $\rho_{w,x,z}$ from $A$ to $D.$
As for intuitive visualization, the equivalence from $D$ to $C$ 
will follow from Figure 2
for the second-fitting step of our assembly algorithm.
For that figure we will use different mappings,
which do not require passing through the block $A.$\\ 

This block equivalence theorem proves an assertion 
we claimed in the introduction: that we are interpreting 
a sum of third powers by assembling 4 copies of a single block.
Each of our blocks $B,C,D,$ in its defined location,
is a copy of the block $A$ in its defined location,
under mappings of rotation and translation.\\

For our assembly algorithm, 
we will also need to label certain axes for blocks 
as \emph{unit cube vacancy axis} or \emph{unit cube overhang axis}. 

\begin{definition}For the blocks $C$ and $D$ only:

The \emph{unit cube vacancy axis} for block $C$ is the axis $y.$

The \emph{unit cube overhang axis} for block $D$ is the axis $z.$
\hfill {\rule{2pt}{1em}}
\end{definition}

These axes help explain our selection of
the particular rotation $\rho_{w,x,z}\thinspace ,$
in our proof of the block equivalence from $A$ to $D.$
The reason for our selection is to ensure: 
that the resulting mapping $\mu_{(D, C)}$ 
which sends $D$ onto $C$ (via $A$) in their defined locations, 
\begin{equation}
\mu_{(D, C)} =  \tau_{x+1, y+1}     \circ 
                \rho_{z, x, y}      \circ 
                (\rho_{w, x, z})^{-1} \circ
                (\tau_{z-1})^{-1}\thinspace ;
\end{equation}
or equivalently, to ensure that 
the resulting mapping $\mu_{(D,C)}^{\prime}$
which sends $D$ onto $C$ (via $A$) in their canonical locations,
\begin{equation}
\mu_{(D, C)}^{\prime} =  \rho_{z, x, y}      \circ 
                        (\rho_{w, x, z})^{-1}\thinspace ; 
\end{equation}
will satisfy properties for the block pair $(D, C)$
which we will require for the second fitting. 
These properties are: 
when $D$ and $C$ are in their canonical locations,
\begin{itemize}
\item 
$\mu_{(D,C)}^{\prime}$ maps the block $D$ 
onto the block $C$ as an isometry;
\item
$\mu_{(D,C)}^{\prime}$ maps the indicator axis $w$ for $D,$ 
to the indicator axis $x$ for $C$;
\item
$\mu_{(D,C)}^{\prime}$ maps the unit cube overhang axis $z$ for $D,$
to unit cube vacancy axis $y$ for $C.$
\end{itemize}
One checks that this mapping
$\rho_{z, x, y} \circ (\rho_{w, x, z})^{-1}$
is the permutation of coordinate axes 
\begin{equation}
\mu_{(D,C)}^{\prime}: (x, y, z, w) \mapsto(w, z, y, x)\thinspace .
\end{equation}         
And this is the only permutation of the axes in $\mathbb{R}^4$
which satisfies that:
it maps $w$ to $x$; and maps $z$ to $y$;
and it is a rotation,  
as a permutation with even parity. 
Thus our selection of the particular rotation $\rho_{w,x,z}$
in the proof of the block equivalence theorem. 

For completeness, we also note an explicit form 
of the mapping $\mu_{(D,C)}\thinspace ,$
\begin{equation}
\mu_{(D,C)}:(x,y,z,w) \mapsto (w+1, z+2, y, x)\thinspace .
\end{equation}

\subsection{The first fitting}

Now we define a rectangular region in $\mathbb{R}^4,$ 
whose (hyper-)volume is $n^2 (n+1)^2.$ We call this region $R_3(n),$ 
because we use it to interpret the sum of third powers. 

First we view the whole space, $\mathbb{R}^4,$ as a union of 
\emph{unit cube sized locations:} 
\begin{equation}
\begin{split}
\mathbb{R}^4 = \bigcup_{a,b,c,d \in \mathbb{Z}}
\{(x,y,z,w) \mid \thickspace
      a-1 \le x \le a;\thickspace
      b-1 \le y \le b;\thickspace
    & c-1 \le z \le c;\thickspace\\
    & d-1 \le w \le d\}.
\end{split}
\end{equation}    
And assign to each of these unit cube sized locations 
a corresponding location label, $(x,y,z,w) = (a,b,c,d),$ 
which records the right endpoint 
of its unit interval on each of the coordinate axes,
just as we did for the unit cubes which satisfy the form (4).

\begin{definition}
We define region $R_3(n) $ 
to be the union of the unit cube sized locations  
whose positions have form as in (9), 
whose location labels $(x,y,z,w)$ satisfy restrictions:
\begin{equation*}\hspace*{4.120em}R_3(n) = 
  \bigcup_{ \substack{x,y \in \{1, 2, \ldots, n+1\}\\
                      z,w \in \{1, 2, \ldots, n  \} } }
(\textrm{unit cube sized location at } (x,y,z,w)).\hspace*{7.940em}
    {\rule{2pt}{1em}} 
\end{equation*}
\end{definition}

The goal of our assembly algorithm is to occupy fully the region 
$R_3(n)$ with the unit cubes which constitute the blocks $A,B,C,D.$

For the first step of the algorithm, the ``first fitting,"
we assemble these blocks in their defined locations.
We will see that the blocks occupy $R_3(n)$ almost entirely, 
but with a layer of overhang unit cubes that lie outside $R_3(n),$ 
and a layer of unit cube sized vacancies that lie within $R_3(n).$

\begin{theorem}[The first fitting]\label{Thm:2.2}
The blocks $A, B, C, D,$ each in its defined location,
do not intersect, except trivially where adjacent.
And the unit cubes of these assembled blocks occupy fully 
the region $R_3(n),$ with the following exceptions.

Lying outside of
$R_3(n)$ is an overhang layer of 
$1^2 + 2^2 + \cdots + n^2$ unit cubes,
at the subcomponent of block $D$ at $z=0$.

And within $R_3(n)$ is a layer of $1^2 + 2^2 + \cdots + n^2$
unit cube sized vacancies at $y=1;$
they lie adjacent to the subcomponent of block $C$ at $y=2$.

\end{theorem}

\begin{proof}
View the whole space $\mathbb{R}^4$ as the union of 
unit cube sized locations stated in (9).   
And assign to each of these unit cube sized locations its corresponding
location label $(x, y, z, w),$ where $x,y,z,w \in \mathbb{Z}.$
 
Each label satisfies exactly one of the following four conditions;
these conditions are mutually exclusive 
and they cover all possible cases,
for the maximum of the set of values $\{x, y, z, w\}$
where any of $x,y,z,w$ may be equal to each other.
\begin{enumerate}
\item[$\alpha$.]
$\max\{x,y,z,w\} = z$.
\item[$\beta$.]
$\max\{x,y,z,w\} = y$,\hspace{0.150em} and\hspace{0.1250em}
$\max\{x,y,z,w\} \neq z, w, x$.
\item[$\gamma$.] 
$\max\{x,y,z,w\} = x$,\hspace{0.050em} and\hspace{0.1250em}
$\max\{x,y,z,w\} \neq z, w$.
\item[$\delta$.]
$\max\{x,y,z,w\} = w$,\hspace{0.06250em} and\hspace{0.1250em}
$\max\{x,y,z,w\} \neq z$.
\end{enumerate}

Now one checks from the
definition of block $A$:
each of the unit cubes that constitute block $A$, 
in its defined location,
has location label that satisfies condition $\alpha$. 
And by a similar check:
each of the unit cubes that constitute block $B$ (respectively $C, D$) 
in its defined location,
has location label that
satisfies condition $\beta$ (respectively $\gamma, \delta$). 

So the blocks $A, B, C, D,$ each in its defined location, 
do not intersect, except trivially where adjacent.
That is, if a unit cube $c_u$ from one block has non-empty intersection 
with a unit cube $c_u^{\prime}$ from a different block:
then they will share a common vertex point
(or unit-interval edge, 
 or unit-square face, 
 or unit 3-cube facet)
through adjacency;
but they will not share the same unit 4-cube location.
\\

Next, the region $R_3(n)$ is the union of 
the unit cube sized locations in (9)
whose location labels $(x,y,z,w)$ satisfy
$x, y \in \{1,2, \ldots, n+1\}$ and
$z, w \in \{1,2, \ldots, n\}.$
We observe the following properties of set containment. 
\begin{itemize}
\item
One checks that the set of all the unit cube sized locations within
the region $R_3(n)$ which satisfy the condition $\gamma,$ is the union
\begin{equation}
\bigcup_{i=1}^n
\bigcup_{\substack{x = i+1\\
                  y \in \{1, 2, \ldots i+1\}\\
                  z,w \in \{1, 2, \ldots, i\}
                 } 
        }
( \text{unit cube sized location at } (x, y, z, w ) ).
\end{equation} 

All the unit cubes that constitute the block $C,$ in its defined location,
lie within this subcomponent of region $R_3(n).$

Conversely: 
all of the unit cube sized locations within the region $R_3(n)$ 
which satisfy condition $\gamma,$
are occupied by unit cubes that define block $C,$ 

-- except not occupied are the unit cube sized locations in (10) at $y=1,$
by the definition of block $C.$

And one checks that for every tuple $(x,y,z,w)$ in (10) where $y=1,$
then there is a corresponding tuple $(x,y,z,w)$ in (10) where $y=2$
with the identical values for $x,z,w;$ and vice versa.
So these unit cube sized vacancies in region $R_3(n)$ at $y=1$ 
are lying adjacent to the unit cubes in 
block $C$ at $y=2,$ and vice versa.
And the number of these unit cube sized vacancies is 
$1^2 + 2^2 + \ldots + n^2.$
\\
 
\item
One checks that the set of all the unit cube sized locations within
the region $R_3(n)$ which satisfy the condition $\delta,$
is the union
\begin{equation}
\bigcup_{i=2}^n
\bigcup_{\substack{w=i\\
                   x,y \in\{1, 2, \ldots, i\}\\
                   z\in\{1, 2, \ldots, i-1\}
                  }
        }
( \text{unit cube sized location at } (x, y, z, w ) ).
\end{equation} 

All the unit cubes that constitute the block $D,$ in its defined 
location, lie within this subcomponent of region $R_3(n),$

-- 
except not lying within $R_3(n)$ are the unit cubes in 
the subcomponent of block $D$ at $z = 0$: 
they form an overhang that lies outside $R_3(n).$ 
More explicitly, these unit cubes are
\begin{equation}
\bigcup_{i=1}^n
\bigcup_{\substack{w=i\\
                   x,y\in\{1, 2, \ldots, i\}\\
                   z=0
                  }
        }
( \text{unit cube at } (x, y, z, w ) ).
\end{equation}
The number of these overhang unit cubes 
is $1^2 + 2^2 + \ldots + n^2.$ 

Conversely: 
all of the unit cube sized locations in the region $R_3(n)$ 
which satisfy condition $\delta,$
are occupied by unit cubes that define block $D.$
\\
 
\item
By similar reasoning, the set of all the unit cube sized locations within 
region $R_3(n)$ which satisfy the condition $\alpha,$ 
is the union
\begin{equation}
\bigcup_{i=1}^n
\bigcup_{\substack{z=i\\
                   x,y,w\in\{1, 2, \ldots i\}
                  }
        }
( \text{unit cube sized location at } (x, y, z, w ) ).
\end{equation} 

All the unit cubes that constitute the block $A$, in its defined location, 
lie within this subcomponent of region $R_3(n).$ 

Conversely: 
all of the unit cube sized locations in the region $R_3(n)$ 
which satisfy condition $\alpha,$
are occupied by unit cubes that define block $A.$
\\

And the set of all the unit cube sized locations within the region $R_3(n)$
which satisfy the condition $\beta,$
is the union
\begin{equation}
\bigcup_{i=1}^n
\bigcup_{\substack{y=i+1\\
                   x,z,w\in\{1, 2, \ldots, i\}
                  }
        }
( \text{unit cube sized location at } (x, y, z, w ) ).
\end{equation}

All the unit cubes that constitute the block $B$, in its defined location, 
lie within this subcomponent of region $R_3(n).$ 

Conversely: 
all of the unit cube sized locations within the region $R_3(n)$ 
which satisfy condition $\beta,$
are occupied by unit cubes that define block $B.$
\end{itemize}
 
So the unit cubes in the blocks $A, B, C, D,$  
each in its defined location, occupy fully the region $R_3(n),$  
with the exception of the layer of $1^2 + 2^2 + \ldots + n^2$
overhang unit cubes in block $D$ at $z=0$
which lie outside $R_3(n),$
and the layer of $1^2 + 2^2 + \ldots + n^2$
unit cube sized vacancies within $R_3(n)$ at $y=1$ 
that lie adjacent to the subcomponent of block $C$ at $y=2.$
\end{proof}

See Figure 1 (at its step 2), for illustration of the first fitting of the blocks.\\

Our rationale for labeling certain axes for blocks 
as unit cube vacancy axis or unit cube overhang axis is now clear: 
from this first fitting theorem, the unit cubes of blocks $A, B, C, D$ 
after the first fitting will occupy the entire region $R_3(n),$ 
with certain exceptions which are stated naturally in that terminology.
\begin{itemize}
\item
 There is a layer of $1^2 + 2^2 + \ldots + n^2$ overhang unit cubes, 
 at the subcomponent of block $D$ where the location label's coordinate value
 equals $0$ on the \emph{unit cube overhang axis} for $D.$
 
\item
 There is a layer of $1^2 + 2^2 + \ldots + n^2$ unit cube sized vacancies, 
 lying adjacent to the subcomponent of block $C$ 
 where the location label's coordinate value equals $2$ 
 on the \emph{unit cube vacancy axis} for $C,$ 
 where the location label's coordinate value on that axis equals $1.$
\end{itemize}

\subsection{The second fitting}

For the second step of our assembly, 
the ``second fitting," we map the 
$1^2 + 2^2 + \ldots + n^2$ overhang unit cubes 
lying outside $R_3(n),$ onto 
the $1^2 + 2^2 + \ldots + n^2$ unit cube sized vacancies 
within $R_3(n).$

First we name certain sets of unit cubes of the blocks $D$ and $C ,$ 
and a certain set of unit cube sized locations.

\begin{definition} For the blocks $D$ and $C$ only,
 each in their defined location:

The \emph{unit cube overhang subcomponent} of $D$ is 
the union of the 
$1^2 + 2^2 + \ldots + n^2$ 
unit cubes in block $D$ 
which have location label coordinate $z=0.$

The \emph{unit cube vacancy-neighbor subcomponent} of $C$ is 
union of 
$1^2 + 2^2 + \ldots + n^2$ 
unit cubes in block $C$
which have location label coordinate $y=2.$

The \emph{unit cube vacancy supercomponent} of $C$ is 
the union of the
$1^2 + 2^2 + \ldots + n^2$
unit cube sized locations in the region $R_3(n)$
which have location label coordinate $y=1 $
and lie adjacent to the unit cube vacancy-neighbor subcomponent of $C.$
\hfill {\rule{2pt}{1em}}

\end{definition}

Recall that
we constructed the rotation $\mu_{(D,C)}^{\prime}$ from $D$ to $C,$
to satisfy certain properties for the block pair $(D, C).$ 
Those properties were: 
when $D$ and $C$ are each in their canonical location,
\begin{itemize}
\item 
$\mu_{(D,C)}^{\prime}$
maps the block $D$ onto the block $C$ as an isometry;
\item
$\mu_{(D,C)}^{\prime}$
maps the indicator axis $w$ for $D,$ 
to the indicator axis $x$ for $C;$
\item
$\mu_{(D,C)}^{\prime}$
maps the unit cube overhang axis $z$ for $D,$ 
to unit cube vacancy axis $y$ for $C.$
\end{itemize}
Now we will see that the rotation $\mu_{(D,C)}^{\prime}$
fulfills another essential role.
When $D$ and $C$ are 
each in their canonical location,
\begin{itemize}
\item 
 $\mu_{(D,C)}^{\prime}$
 maps the unit cube overhang subcomponent of $D$
 onto the unit cube vacancy-neighbor subcomponent of $C,$ as a
 bijection of unit cubes. 
\end{itemize}
We state and prove this as a theorem,
more explicitly in the context 
where $D$ and $C$ are each in their defined location.
Recall the formula for the mapping $\mu_{(D,C)}\thinspace ,$
\begin{equation*}
\mu_{(D,C)} = \tau_{x+1, y+1}\circ
              \rho_{z,x,y}\circ
              (\rho_{w,x,z})^{-1}\circ
              (\tau_{z-1})^{-1}.
\end{equation*}
\begin{theorem}[Replication]\label{Thm:2.3}
When the blocks $D,C$ are each in their defined location,
the mapping $\mu_{(D,C)}$
sends the unit cube overhang subcomponent of $D$
onto the unit cube vacancy-neighbor subcomponent of $C,$
as a bijection of unit cubes.
\end{theorem}
\begin{proof}
The unit cube overhang subcomponent of $D$ is as stated in (12):
\begin{equation*}
\bigcup_{i=1}^n 
\bigcup_{\substack{w=i\\
                   x,y\in\{1, 2, \ldots, i\}\\
                   z=0
                  }
        }
                 ( \textrm{unit cube at } (x,y,z,w) )\thinspace .                   
\end{equation*} 
This is the subcomponent of block $D,$ 
in its defined location,
at its smallest location label coordinate for $z$: at $z=0.$

Now applying the translation $\tau_{z-1}^{-1}$ yields
\begin{equation*}
\bigcup_{i=1}^n 
\bigcup_{\substack{w=i\\
                   x,y\in\{1, 2, \ldots, i\}\\
                   z=1
                  }
        }
                 ( \textrm{unit cube at } (x,y,z,w) )\thinspace .                   
\end{equation*}
This is the subcomponent of block $D,$ 
in its canonical location,
at its smallest location label coordinate for $z$: at $z=1.$

Then applying the rotation $\rho_{z,x,y} \circ (\rho_{w,x,z})^{-1}$ yields
\begin{equation*}
\bigcup_{i=1}^n 
\bigcup_{\substack{x=i\\
                   w,z\in\{1, 2, \ldots, i\}\\
                   y=1
                  }
        }
                 ( \textrm{unit cube at } (x,y,z,w) )\thinspace .                   
\end{equation*}
This is the subcomponent of block $C,$ 
in its canonical location,
at its smallest location label coordinate for $y$: at $y=1.$

Then applying the translation $\tau_{x+1,y+1}$ yields
\begin{equation*}
\bigcup_{i=1}^n 
\bigcup_{\substack{x=i+1\\
                   w,z\in\{1, 2, \ldots, i\}\\
                   y=2
                  }
        }
                 ( \textrm{unit cube at } (x,y,z,w) )\thinspace .                   
\end{equation*}
This is the subcomponent of block $C,$ 
in its defined location,
at its smallest location label coordinate for $y$: at $y=2.$
This final union of unit cubes is the 
unit cube vacancy-neighbor subcomponent of $C.$ 
\end{proof}

When $D$ and $C$ are each in their defined location,
we say that this action of $\mu_{(D,C)}$ on the unit cube overhang
subcomponent of $D$ 
serves to \emph{replicate} the unit cube
vacancy-neighbor subcomponent of $C.$
Hence we call this theorem ``Replication."

In other words, 
when $D$ and $C$ are each in their canonical location,
the rotation $\mu_{(D,C)}^{\prime}$ 
sends the subcomponent of $D$ at $z=1$
onto the subcomponent of $C$ at $y=1.$ 
\\

We could have proved this replication theorem directly, 
but with less visual insight, 
by applying in one step the mapping
\hspace*{0.1250em}
$\mu_{(D,C)}: (x, y, z, w) \mapsto (w+1, z+2, y, x),$
to the unit cubes that define the overhang subcomponent of $D.$
\\

Now define a mapping $\phi_{C, u}\thinspace ,$ 
which serves to \emph{inflate} outward 
the replicated unit cube vacancy-neighbor subcomponent of $C,$ 
onto the unit cube vacancy supercomponent of $C.$

\begin{definition} For the block $C$ only,
in its defined location:\\
We define the \emph{unit cube inflation map} $\phi_{C, u}\thinspace ,$
to be the function that acts on the unit cube vacancy-neighbor
subcomponent of block $C,$ 
by translating each unit cube by $-1$ 
in the unit cube vacancy axis of $C,$ i.e. from $y=2$ to $y=1.$
\hfill {\rule{2pt}{1em}}
\end{definition}

When we apply the inflation map $\phi_{C,u}$ to the
replicated unit cube vacancy-neighbor subcomponent of $C ,$
it sends these unit cubes in bijection onto
the unit cube sized vacancies in the 
unit cube vacancy supercomponent of $C.$

We record this final step of our assembly algorithm as a theorem.
\begin{theorem}[The second fitting]\label{Thm:2.4}
When the blocks $D,C$ are each in defined location,
the mapping $\phi_{C,u} \circ \mu_{(D,C)}$
of replication and inflation
sends the unit cube overhang subcomponent of $D$ onto the
unit cube vacancy supercomponent of $C ,$ 
as a bijection of unit cubes onto unit cube sized locations.
\end{theorem}

See Figure 1 (at its step 3a) 
and Figure 2 (at its step 3b) 
for illustrations of this 
second-fitting mapping
\begin{equation*}
\phi_{C,u}\circ\mu_{(D,C)}:(x,y,z,w)\mapsto(w+1, z+1, y, x)
\end{equation*}
from the unit cube overhang subcomponent of the block $D$ 
onto the unit cube vacancy supercomponent of the block $C.$\\

For the illustrations in those figures, we 
use a composition of different rotations and translations
to achieve this second-fitting mapping. Namely
$(\rho_2 \circ \tau)\circ(\rho_1),$ where
\begin{alignat*}{4}
\rho_1 &: (x, y, z, w) \mapsto (w   &&,\medspace n+1-y&&,\medspace z    &&,\medspace x)\\
\tau   &: (x, y, z, w) \mapsto (x+1 &&,\medspace y    &&,\medspace z    &&,\medspace w)\\
\rho_2 &: (x, y, z, w) \mapsto (x   &&,\medspace z+1  &&,\medspace n+1-y&&,\medspace w).
\end{alignat*} 
One checks that 
\begin{equation*}
\rho_2 \circ \tau \circ \rho_1 : (x, y, z, w) \mapsto
                                 (w+1, z+1  , y, x),
\end{equation*}
so this equals $\phi_{C,u} \circ \mu_{(D,C)}$ as a mapping 
from $\mathbb{R}^4$ to $\mathbb{R}^4$.

One reason for selecting these new rotations and translation 
is they make easier to visualize
the mapping from the unit cube overhang subcomponent of $D$
onto the unit cube vacancy supercomponent of $C$: 
in Figure 2 we illustrate the individual mappings 
$\rho_1, \tau, \rho_2$
each as (one or a sequence of) isometry transformations 
in a plane parallel to 2 of our coordinate axes.

Another reason for selecting these new rotations and translation 
is they make easier to visualize 
the mapping $\mu_{(D,C)},$
\begin{equation*}
\mu_{(D,C)} = (\phi_{C,u})^{-1} \circ (\rho_2 \circ \tau \circ \rho_1),
\end{equation*}
which sends the subcomponent of $D$ at $z=0$ 
onto the subcomponent of $C$ at $y=2$ 
(when each block is in its defined location).

And from this, one checks that more generally:  
$\mu_{(D,C)}$ 
sends the subcomponent of $D$ at $z = i-1,$
onto the subcomponent of $C$ at  $y = i+1,$ 
for every $i \in \{1, 2, \ldots, n\}$ 
(when each block is in its defined location).
This completes a visualization that 
the whole blocks $D$ and $C$ are equivalent under
mappings of rotations and translations.
\\

This completes our filling of the rectangular region:
all of the unit cubes of the blocks $A,B,C,D$ map bijectively,
through the first fitting and the second fitting,
onto the unit cube sized locations in $R_3(n).$
And this concludes our interpretation of the identity that
$\thinspace (1^3 + 2^3 + \cdots + n^3)\cdot 4 =\thinspace n^2 (n+1)^2\thinspace .$\\

See Figure 2 (at its step 4) for illustration of the filled region $R_3(n).$


\subsection{A q-analogue of the sum of third powers, 
            from the four-block assembly}

As a corollary we present a $q$-analogue 
of the sum of third powers equation (2), 
in the algebraically equivalent form
$(1^3 + 2^3 + \ldots + n^3)\cdot 4 = n^2 (n+1)^2.$

Notation: for any non-negative integer $k,$
we write the $q$-analogue of $k$ as
\begin{equation*}
[k]_q\thinspace =\thinspace {(1-q^k)} / {(1-q)}\thinspace 
      =\thinspace 1 + q + q^2 + \ldots + q^{k-1}.
\end{equation*}
We are not necessarily interested in 
evaluation at particular values of the parameter $q,$ but note that
$\lim_{q \to 1^{-}} [k]_q = k.$
We also write $q$-analogues of 
the factorial and the binomial coefficient:
\begin{equation*}
[k]_q ! \thinspace
= \thinspace
\prod_{i=1}^{k} [i]_q\medspace ,
\hspace*{1.50em} [0]_q! \thinspace
= \thinspace
1\medspace ,
\hspace*{2.00em}
{\genfrac{[}{]}{0pt}{0}{n}{k}}_q\thinspace
=\thinspace 
\frac{[n]_q !}{[k]_q ! \medspace [n-k]_q !}\medspace .
\end{equation*}

Begin by recalling our four-block assembly
for interpreting the sum of third powers, 
from earlier in this section. 
Now for each unit cube $c_u$
in the final fitted configuration in the region $R_3(n)$,
take its location label $(x,y,z,w) = (a,b,c,d)$
and assign the monomial
\begin{equation*}
q^{(a-1) + (b-1) + (c-1) + (d-1)}.
\end{equation*}
Here the exponent of $q$ is the \emph{taxicab distance} 
between $c_u$ and the unit cube in $R_3(n)$
whose location label is $(1,1,1,1).$
That is, this exponent is the sum of the distances along each 
of the individual coordinate axes,
between $(a,b,c,d)$ and $(1,1,1,1).$
\\

Summing this monomial over all unit cubes
in this fitted assembly in $R_3(n)$ gives
\begin{equation*}
(q^0 + q^1 + \ldots + q^{n})\thinspace
(q^0 + q^1 + \ldots + q^{n})\thinspace
(q^0 + q^1 + \ldots + q^{n-1})\thinspace
(q^0 + q^1 + \ldots + q^{n-1}),
\end{equation*}
since $a$ and $b$ take every value in $\{1, 2, \ldots, n+1\},$
and $c$ and $d$ take every value in $\{1, 2, \ldots, n\},$
for this rectangular region $R_3(n) .$
One checks that this expression equals 
$[n+1]_q\thinspace
 [n+1]_q\thinspace
 [n]_q\thinspace
 [n]_q.$\\

Now instead, sum this monomial
over the unit cubes in this assembly 
that originated in block $A,$
and likewise for $B, C,$ and $D$ 
(whose layer at $z=0$ got relocated 
by the second-fitting mapping). 
We obtain for the sum over the unit cubes in block $A$
\begin{equation*}
\sum_{i=1}^{n}\thinspace q^0[i]_q\cdot
               q^0[i]_q\cdot
               q^{i-1}\cdot
               q^0[i]_q ,
\end{equation*}
which follows from the definition of this block.
That is: in the definition of block $A,$ 
when the union's index $i$ takes a given value, 
then for each of the $i^3$ unit cubes $c_u$ in $A$
which correspond to that $i,$
the $c_u$ has a location label $(a,b,c,d)$ 
where $c$ takes value $i,$ 
and $a,b,d$ each take any value in $\{1, 2, \ldots, i\}.$
So for that $i,$ 
the sum of the corresponding values of the monomial
$q^{(a-1) + (b-1) + (c-1) + (d-1)}$ is equal to
\begin{equation*}
(q^0 + q^1 + \ldots + q^{i-1}) \cdot
(q^0 + q^1 + \ldots + q^{i-1}) \cdot
(q^{i-1}) \cdot
(q^0 + q^1 + \ldots + q^{i-1}),
\end{equation*}
which equals the summand 
$\thinspace q^0[i]_q \cdot q^0[i]_q \cdot q^{i-1} \cdot q^0[i]_q\thinspace $
of our displayed summation.

By similar reasoning,
we obtain for the sum over all of the fitted blocks $A, B, C, D$:
\begin{alignat*}{2}
\sum_{i=1}^{n}\thinspace \Big(&\big(q^0[i]_q\thinspace
               q^0[i]_q\thinspace
               q^{i-1}\thinspace
               q^0[i]_q\big) 
               +
               \big(q^0[i]_q\thinspace
               q^{i}\thinspace
               q^0[i]_q\thinspace
               q^0[i]_q\big) 
               +
               \big(q^{i}\thinspace
               q^1[i]_q\thinspace
               q^0[i]_q\thinspace
               q^0[i]_q\big)\\
               &+\big(
               q^0[i]_q\thinspace
               q^0[i]_q\thinspace
               q^0[i-1]_q\thinspace
               q^{i-1} 
               +
               q^{i}\thinspace
               q^0\thinspace
               q^0[i]_q\thinspace
               q^0[i]_q \big)\Big).
\end{alignat*}
Here we have applied that for the block $D ,$ 
its component at $z=0$ gets mapped in the second fitting,
onto the set of unit cube sized locations at $y=1$ that 
lie adjacent to the subcomponent of block $C$ at $y=2 .$

Now by equating,
our $q$-analogue of the sum of third powers 
formula is established:
\begin{theorem}\label{Thm:2.5}
\begin{equation}
\sum_{i=1}^{n}\thinspace
(q^{i-1} + q^i + q^{i+1})\thinspace
[i]_q^3\thinspace
+\thinspace
\sum_{i=1}^{n}\thinspace
(q^{i-1}[i-1]_q + q^i \cdot 1)\thinspace
[i]_q^2\thinspace
=\thinspace
[n+1]_q^2\thinspace 
[n]_q^2.
\end{equation}
\end{theorem}
One checks that, by taking $\lim_{q\to 1^{-}}$ 
of each side of this equation of polynomials, 
we will then recover the integer equation
$(1^3 + 2^3 + \ldots + n^3) \cdot 4 = n^2 (n+1)^2 .$


\section{The one-block construction}
In this section we interpret the identity (2) in its initial form,
\begin{equation*}
1^3 + 2^3 + \ldots + n^3 = \tfrac{1}{4} n^2 (n+1)^2 ,
\end{equation*}
by applying a partition and a mapping reassembly 
to a single 4-dimensional block.

Begin by recalling from Section 2, the definitions
for the unit cube $c_u\thinspace ,$
the location label $(x,y,z,w) = (a,b,c,d) ,$
the blocks $A, B, C, D, $ 
and the rectangular region $R_3(n) .$
\\

Our construction employs 4 steps. 
We state them here, and give proofs afterward.
\begin{enumerate}
\item 
Select the block $B .$
It comprises $\sum_{i=1}^{n} i^3$ unit cubes.
\item
Partition $B$ into two parts:
$B_{c\ge d}$ which comprises the unit cubes $c_u$ in $B$
whose location labels satisfy $c\ge d ;$ and
$B_{c<d}$ comprising those that satisfy $c<d .$
\item
For the part $B_{c<d}\thinspace ,$
apply the isometry mapping
\begin{equation}
\psi_{z,w,y-1,x} : (x,y,z,w) \mapsto (z,w,y-1,x)\thinspace .
\end{equation}
We write:
$\thinspace \psi_{z,w,y-1,x}(B_{c<d}) = A_{a<b}\thinspace .$
[See notes below.]
\item 
Assemble the disjoint parts $B_{c\ge d}$ and $A_{a<b}\thinspace .$
Their union is the set of unit cubes in 
region $R_3(n)$ whose location labels satisfy 
$1\le a < b \le n+1$ and $1 \le d \le c \le n .$
We write: $B_{c\ge d} \bigcup A_{a<b} = S_{3,3}(n).$
This comprises ${\binom{n+1}{2}}^2$ unit cubes. 
[See notes below.]
\\
\end{enumerate}

\noindent
Notes for Step 3 of the one-block construction:

\begin{lemma}\label{Thm:3.1}
The mapping $\psi_{z,w,y-1,x}\thinspace ,$ as defined in (16),
is an isometry which sends 
the block $B$ onto block $A$ 
as a bijection of unit cubes. In particular, it sends 
the unit cubes in $B$ whose location labels satisfy $c < d ,$ onto 
the unit cubes in $A$ which satisfy $a < b .$
\end{lemma}
\begin{proof}
This $\psi_{z,w,y-1,x}$ is an isometry in $\mathbb{R}^4 $
because it is equal to 
$\phi_{w,y} \circ \rho_{z,x,y} \circ \tau_{y-1}\thinspace ,$ where
$\tau_{y-1}$ is the translation by $-1$ in $y ,$ and
$\rho_{z,x,y} \thinspace : (x,y,z,w) \mapsto (z,x,y,w) $ 
is a rotation, 
and
$\phi_{w,y} \thinspace : (x,y,z,w) \mapsto (x,w,z,y) $ is a reflection.

And this translation $\tau_{y-1}$ maps the block $B$
from its defined location to canonical location.
Then rotation $\rho_{z,x,y}$ maps $B$ onto $A ,$
for each block in its canonical location.
Then $\phi_{w,y}$ sends $A$ to itself.
This canonical location for $A$ is equal to its defined location.

In particular, $\psi_{z,w,y-1,x}$ as defined in (16)
maps the unit cubes in $B$ whose location labels satisfy $z<w\thinspace $ 
(that is, the unit cubes in $B_{c < d\thinspace }$),
onto the unit cubes in $A$ which satisfy $x < y .$
By our location labelling, these are the unit cubes in $A$ where $a<b .$ 
\end{proof}

\noindent
Notes for Step 4 of the one-block construction:
\begin{lemma}\label{Thm:3.2}
The parts $B_{c\ge d}$ and $A_{a<b}$ are mutually disjoint as sets of unit cubes.
And their union is the set of unit cubes $c_u$ in $R_3(n)$ 
whose location labels satisfy both 
$1 \le a < b \le n+1$ and $1 \le d \le c \le n .$
The number of unit cubes in this union is ${\tbinom{n+1}{2}}^2 .$
\end{lemma}

\begin{proof}
These $B_{c \ge d}$ and $A_{a<b}$ are mutually disjoint as sets of unit cubes,
because they are subsets of mutually disjoint sets of unit cubes: 
the blocks $B$ and $A .$

As for their union,
we will structure our argument in terms of set containments.
Recall from the proof of the first fitting Theorem~\ref{Thm:2.2}: 
for each unit cube sized location in region $R_3(n) ,$ 
its location label must satisfy exactly one
of the four conditions $\alpha, \beta, \gamma, \delta .$

For those locations in $R_3(n)$ which satisfy the inequalities 
$a<b$ and $d\le c ,$ they cannot satisfy condition $\gamma$
(which would require $a \ge b$), nor can they satisfy condition $\delta$
(which would require $d > c$). So each such location satisfies either
condition $\alpha$ or $\beta .$

Also recall from the proof of Theorem~\ref{Thm:2.2}: 
\begin{itemize}
\item
The unit cube sized locations in $R_3(n)$ which satisfy condition $\alpha ,$
are all occupied by the unit cubes in the block $A$ in its defined location.
And conversely, all the unit cubes in block $A ,$ in its defined location,
lie within $R_3(n)$ and satisfy condition $\alpha .$
\item
The unit cube sized locations in $R_3(n)$ which satisfy condition $\beta ,$
are all occupied by the unit cubes in the block $B$ in its defined location.
And conversely, all the unit cubes in block $B ,$ in its defined location,
lie within $R_3(n)$ and satisfy condition $\beta .$
\end{itemize}

So for each unit cube sized location in region
$R_3(n)$ whose location label 
satisfies both inequalities $a<b$ and $d \le c ,$ we have that:
\begin{itemize}
\item 
this location is occupied by a unit cube in $A$ 
(which satisfies $d \le c$ by definition),
in the subcomponent of $A$ where $a<b ;$ or
\item 
this location is occupied by a unit cube in $B$ 
(which satisfies $a < b$ by definition),
in the subcomponent of $B$ where $d \le c ;$ 
\item 
but not both -- the blocks $A$ and $B$ are disjoint as sets of unit cubes.
\end{itemize}
So we have a set containment: each of the unit cube sized locations
in $R_3(n)$ whose location label satisfies both $a<b$ and $d\le c ,$
is occupied by a unit cube in the union of the disjoint sets
$A_{a<b}$ and $B_{c \ge d}\thinspace .$

And conversely, we have that:
\begin{itemize}
\item
each of the unit cubes $c_u$ in $A_{a<b}\thinspace ,$
in the defined location of $A ,$ lies within the region 
$R_3(n)$ in its subcomponent where both $a<b$ and $d\le c;$ and
\item 
each of the unit cubes $c_u$ in $B_{c \ge d}\thinspace ,$
in the defined location of $B ,$ lies within the region 
$R_3(n)$ in its subcomponent where both $a<b$ and $d\le c .$ 
\end{itemize}
So the set of location labels of the unit cube sized locations in
$R_3(n)$ which satisfy both $a<b$ and $c\ge d ,$
is equal to the set of location labels of the unit cubes in the union
of the disjoint sets $A_{a<b}$ and $B_{c\ge d} .$
This completes our proof of the set equality,
between the union $B_{c\ge d} \bigcup A_{a<b} \thinspace $
and the set of unit cubes in $R_3(n)$ which satisfy both
$1 \le a < b \le n+1$ and $1 \le d \le c \le n .$
\\

This particular union of unit cubes in $R_3(n)$
has the shape of a step triangle in the
$(x,y)$-plane, for each value of $(z,w) .$
And its shape is a step triangle in the $(z,w)$-plane,
for each value of $(x,y) .$

In the next section (in Lemma~\ref{Thm:3.4}), we will show that
the number of unit squares in the $(z,w)$-plane whose location labels
$(c,d)$ satisfy the inequality $1 \le d \le c \le n $
is equal to $\tbinom{n+1}{2} ,$ 
by proving a stronger result: a $q$-analogue of this equation.
And in Lemma~\ref{Thm:3.5} 
we will prove a $q$-analogue of the equation that
the number of unit squares in the $(x,y)$-plane whose location labels
$(a,b)$ satisfy $1 \le a < b \le n+1$ is $\tbinom{n+1}{2} .$ 
\end{proof}

For nomenclature: in $4$-dimensional geometry,
a Cartesian product of two triangles is a 
``triangular -- triangular (3,3) duoprism."
We will call the shape of our assembled solid 
an \emph{(oblique) step triangular -- (oblique) step triangular duoprism}.

And we name 
this set of the unit cubes in $R_3(n)$
which satisfy both $1 \le a < b \le n+1$ and
$1 \le d \le c \le n ,$
to be
$S_{3,3}(n) .$
\\

See Figure 3
for illustration of this one-block construction. 
\\

\begin{figure} 
    \centering

    \includegraphics[scale=0.99800]{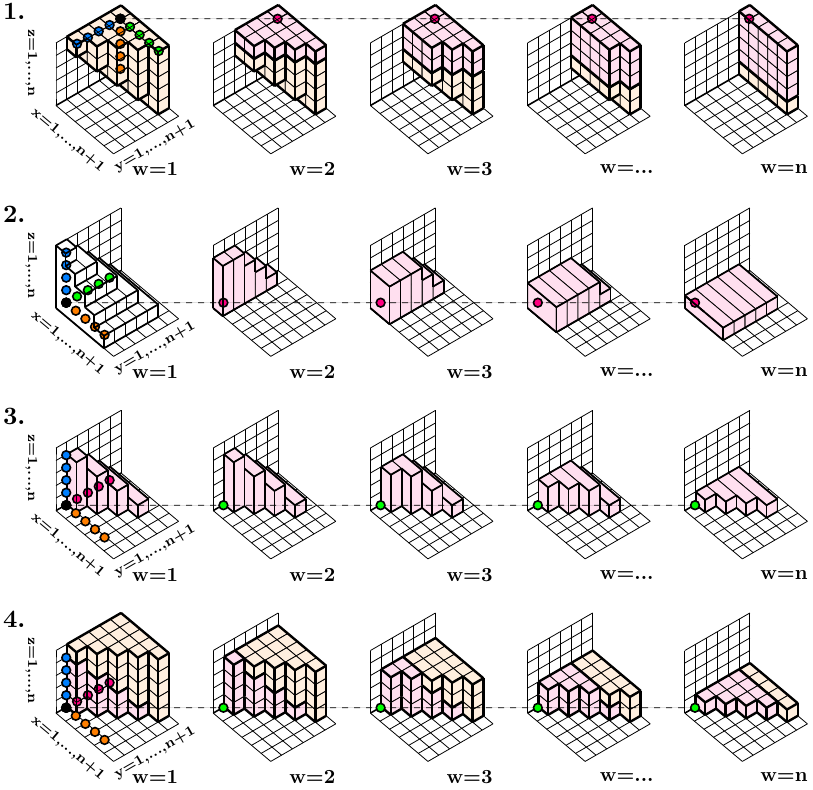}
    \caption{
    We interpret the identity $1^3 + 2^3 + \ldots + n^3 = \tfrac{1}{4} n^2 (n+1)^2 ,$
    using a single $4$-dimensional block.
    We select the block $B ,$ here depicted for $n=5 .$
    This block is a union of $\sum_{i=1}^{n} i^3$ unit hyper-cubes 
    (we call \emph{unit cubes}).
    Denote the location label of an arbitrary unit cube $c_u$ by
    $(x,y,z,w) = (a,b,c,d) .$
    \\
    (1.) Partition $B$ into two parts:
    $B_{c\ge d}$ comprising the unit cubes where $c \ge d$
    (in this figure, hued in brown);
    and $B_{c < d}$ comprising those where $c < d$ (hued in red).
    \\
    (2.) Select $B_{c < d} .$
    Translate it by $\tau_{y-1} : y \mapsto y-1 ;$
    and rotate by $\rho_{z,x,y} : (x,y,z) \mapsto (z, x, y) .$
    The image is the subcomponent of block $A$ whose unit cubes satisfy $a < d .$
    (To assist the visualization, we include an outline of $A$ where $w=1 .$)
    \\
    (3.) Then reflect by $\phi_{w,y} : (y,w) \mapsto (w,y) .$
    The image is the subcomponent of $A$ whose unit cubes satisfy $a < b ;$
    we call it $A_{a < b} .$
    (To assist the visualization, in steps 2 and 3 we have depicted each set
    of unit cubes which share the same coordinate pair in $(y,w) ,$
    as joined into a single subcomponent.)
    \\
    (4.) These $B_{c \ge d}$ and $A_{a < b}$ are mutually disjoint.
    Their union, which we call $S_{3,3}(n) , $ 
    is the subcomponent of the region 
    $R_{3}(n)$ whose unit cubes $c_u$ satisfy that both
    $1 \le a < b \le n+1$ and $1 \le d \le c \le n .$
    The volume of the assembled solid is ${\binom{n+1}{2}}^2 .$
    }
\end{figure}

\begin{remark}
Notes on prior art:

A. T. Benjamin and M. E. Orrison \cite{example97_benjamin_orrison} 
have presented 
combinatorial proofs of the sum of third powers formula,
$\sum_{k=1}^{n}k^3 = {\tbinom{n+1}{2}}^2 .$ 
We summarize their simpler correspondence. 
\begin{itemize}
\item 
They let $S = \{(h,i,j,k) \mid 1 \le h, i, j \le k \le n\} .$
This set has cardinality $\sum_{k=1}^{n} k^3 .$
\item 
They let
$T = \{ ((x_1, x_2), (x_3, x_4)) \mid 1 \le x_1 \le x_2 \le n;\thickspace 
                                      1 \le x_3 \le x_4 \le n\} .$
This set has cardinality ${\tbinom{n+1}{2}}^2 .$
\item 
They define a bijection $g : S \rightarrow T$ by the rule
\begin{equation*}
    g((h, i, j, k)) = \begin{cases} 
          ((h,i),(j,k)) & \textrm{ if } h \le i \\
          ((j,k),(i,h-1)) & \textrm{ if } h > i \thinspace . 
       \end{cases}
\end{equation*}
\end{itemize}

Observe that for our one-block combinatorial proof 
of the sum of third powers formula,
we could have 
achieved the same outcome as we did with our original
selections of sets of unit cubes,
by instead selecting each of their images under a translation by
$-1$ in $y .$
Let $\tau_{y-1} : (x,y,z,w) \mapsto (x,y-1,z,w) .$
Our new construction would be worded as follows.
\begin{itemize}
\item
Select the block $\tau_{y-1}(B) .$
It comprises $\sum_{i=1}^n i^3$ unit cubes.
\item 
Partition $\tau_{y-1}(B)$ into two parts:
$\tau_{y-1}(B_{c \ge d})$ which comprises the unit cubes $c_u$ in 
$\tau_{y-1}(B)$
whose location labels satisfy $c \ge d ;$
and $\tau_{y-1}(B_{c < d})$ comprising those that satisfy $c < d .$

\item
For the part $\tau_{y-1}(B_{c < d})\thinspace ,$
apply the isometry 
$\psi_{z,w-1,y,x} : (x,y,z,w) \mapsto (z, w-1, y, x )\thinspace ,$      
i.e. the composite map 
$\tau_{y-1} \circ \psi_{z,w,y-1,x} \circ \tau_{y-1}^{-1}\thinspace .$

We write: 
$\psi_{z,w-1,y,x}( \tau_{y-1}(B_{c < d}) ) 
= A_{a \le b \le n-1}\thinspace .$
This is the image of $A_{a<b}$ under $\tau_{y-1}\thinspace .$
\item 
Assemble the disjoint parts 
$\tau_{y-1}(B_{c \ge d})$ and 
$A_{a \le b \le n-1}\thinspace .$
Their union is the set of unit cubes in region $R_3(n)$
whose location labels satisfy $1 \le a \le b\le n$ and
$1 \le d \le c \le n .$
That is, their union is the image of our assembled solid 
$S_{3,3}(n)$ under $\tau_{y-1}\thinspace .$
This comprises ${\tbinom{n+1}{2}}^2$ unit cubes.
\end{itemize}

We show that 
these 
correspondences between 
$\sum_{k=1}^n k^3$ and ${\tbinom{n+1}{2}}^2$
-- that of Benjamin and Orrison, and 
   that of our ``$\tau_{y-1}$-translated" one-block construction --
are equivalent.

\begin{enumerate}
\item 
For each 
unit cube $c_u$ in $\tau_{y-1}(B)$
with location label $(a,b,c,d),$
we assign a 
4-tuple $(h,i,j,k)\in S ,$
by the assignment rule
$f_1 : ( a,b,c,d ) \mapsto (d, c, a, b) .$
\item 
For each 
$((x_1, x_2), (x_3, x_4)) \in T ,$
we assign 
a unit cube in $\tau_{y-1}(S_{3,3}(n))$
with 
location label $(x,y,z,w) ,$ by 
the assignment rule 
$f_2 : ((x_1, x_2), (x_3, x_4)) \mapsto (x_3, x_4, x_2, x_1) .$

\item 
Now for each unit cube $c_u$ in $\tau_{y-1}(B)$
with location label $(a,b,c,d) ,$
the image of that unit cube under the composite mapping 
$ f_2 \circ g \circ f_1$ is evaluated by: 

\begin{equation*}
\hspace*{-8.6250em}
f_1 :  (a,b,c,d) \mapsto (d,c,a,b) \thinspace .
\end{equation*}

\begin{equation*}
    g : (d, c, a, b) \mapsto \begin{cases} 
          ((d,c),(a,b)) & \textrm{ if } d \le c \\
          ((a,b),(c,d-1)) & \textrm{ if } d > c \thinspace . 
       \end{cases}
\end{equation*}

\begin{equation*}
\hspace*{-3.6250em}
     \begin{cases} 
          f_2 : ((d,c),(a,b  )) & \mapsto (a, b  , c, d) \\
          f_2 : ((a,b),(c,d-1)) & \mapsto (c, d-1, b, a) \thinspace . 
       \end{cases}
\end{equation*}

One sees that this $f_2 \circ g \circ f_1$
is the same mapping as 
is applied
in 
our ``$\tau_{y-1}$-translated" 
one-block construction:
it fixes the unit cubes 
$c_u$ 
which lie 
in $\tau_{y-1}(B_{c \ge d})\thinspace ,$
and applies the isometry 
$\psi_{z,w-1,y,x}$ to the unit cubes 
$c_u$ 
which lie
in $\tau_{y-1}(B_{c < d}) \thinspace .$
\end{enumerate}

We conclude this remark by 
praising the elegance of the construction of Benjamin and Orrison.
We suggest that our rediscovery of this correspondence is of 
independent interest, for the applications which are 
afforded by its setting in 4-dimensional geometry.
\\
\end{remark}

For corollaries,
now we will interpret the $q$-analogue equations for (2) in its initial form,
$1^3 + 2^3 + \ldots +n^3 = \tfrac{1}{4}n^2 (n+1)^2 ,$
from each of the four sets of authors mentioned in the introduction.

In each instance we will: 
\begin{itemize}
\item
Start with our one-block construction,
which selects the block $B,$
then partitions it into parts $B_{c\ge d}$ and $B_{c<d}\thinspace ,$
then maps $B_{c<d}$ onto $A_{a<b}$ by the isometry
$\psi_{z,w,y-1,x}\thinspace ,$
and then constructs the union 
$B_{c\ge d}\thinspace \bigcup A_{a<b} = S_{3,3}(n)\thinspace $
as an assembled solid in $R_3(n) .$
\item
Select a monomial of the form 
$q\wedge(c_u\textrm{'s taxicab distance to }c_1) ,$
for some point $c_1 ,$ 
which we will sum over 
the unit cubes $c_u$ in $B_{c\ge d}\thinspace .$
\\
Select a monomial of the form 
$q\wedge(c_u\textrm{'s taxicab distance to }c_2) ,$
for some point $c_2 ,$ 
which we will sum over 
the unit cubes $c_u$ in $B_{c < d}\thinspace .$
\\
Select a monomial of the form 
$q\wedge(c_u\textrm{'s taxicab distance to }c_3) ,$
for some point $c_3 ,$ 
which we will sum over 
the unit cubes $c_u$ in $A_{a < b}\thinspace .$

We require here that 
the summation over $B_{c\ge d}$ and
the summation over $A_{a <  b}$ must each be adding the same common monomial.

And we require that 
the summation over $A_{a < b}$ and
the summation over $B_{c < d}$ 
must be equal as sums of their respective monomials.

\item
Interpret the one side of the authors' equation 
which contains their $q$-analogue of ${\binom{n+1}{2}}^2 ,$
as a sum over the unit cubes $c_u$ in the assembled solid
$S_{3,3}(n)\thinspace .$
That is, as the total of 
the summation over the unit cubes $c_u$ in $B_{c \ge d}\thinspace ,$ plus
the summation over the unit cubes $c_u$ in $A_{a < b}\thinspace ,$
of their common monomial. 
 
\item
Interpret the other side of the authors' equation,
which contains their $q$-analogue of $\sum_{i=1}^{n} i^3 ,$
as a sum over the unit cubes $c_u$ in the block $B .$
That is, as the total of 
the summation over the unit cubes $c_u$ in $B_{c \ge d}\thinspace ,$ plus
the summation over the unit cubes $c_u$ in $B_{c < d}\thinspace ,$ of their 
respective monomials.
\end{itemize}


\subsection{Interpreting the q-analogue
         from Garrett and Hummel}

In this section we interpret the $q$-analogue of 
the sum of third powers equation from 
Garrett and Hummel~\cite{example10_garrett_etal}. 
Garrett and Hummel's equation states:
\begin{equation}
\sum_{k=1}^{n} \thinspace
               q^{k-1} 
               \left(\frac{1-q^k}{1-q}\right)^2
               \left[ \frac{1-q^{k-1}}{1-q^2} + 
                      \frac{1-q^{k+1}}{1-q^2} \right]\thinspace
               =\thinspace
               {\genfrac{[}{]}{0pt}{0}{n+1}{2}}_q^{2}\thinspace .
\end{equation}
To enable interpretation with taxicab distance, 
we re-express the left side with algebra.
\begin{equation}
\sum_{k=1}^{n}  \thinspace
                [k]_q\thinspace
                q^{k-1}
               \frac{[k+1]_q\thinspace [k]_q}{1+q} \thinspace
+ \thinspace
\sum_{k=1}^{n}  \thinspace
                [k]_q\thinspace
                q^{k-1}
               \frac{[k]_q\thinspace [k-1]_q}{1+q}\thinspace
               = \thinspace
               {\genfrac{[}{]}{0pt}{0}{n+1}{2}}_q^{2}\thinspace .
\end{equation}

We begin our interpretation of (18) by selecting the block $B .$
Let $c_u$ denote an arbitrary unit cube in $B ,$ 
and let its location label be $(x,y,z,w) = (a,b,c,d).$
Partition $B$ into two parts, $B_{c \ge d}$ and $B_{c < d}\thinspace ,$ 
as defined at the start of this section 3.

For this $q$-analogue equation, 
we will select for our $q$-counts of taxicab distances: 
\begin{itemize}
\item
$\hspace*{0.250em}\hphantom{ q^{-1} \cdot }\sum_{c_u\in B;\thinspace c \ge d}\thinspace 
     q^{(a-1) + (b-2) + (c-1) + (d-1)}$
\item
$q^{-1} \cdot \sum_{c_u\in B;\thinspace c < d}\thinspace 
     q^{(a-1) + (b-2) + (c-1) + (d-1)}$
\item
$q^{-1} \cdot  \sum_{c_u\in A; a < b}\thinspace 
     q^{(a-1) + (b-1) + (c-1) + (d-1)} \thinspace
 = \thinspace
 \sum_{c_u\in A;\thinspace a < b}\thinspace 
     q^{(a-1) + (b-2) + (c-1) + (d-1)} \thinspace .$     
\end{itemize}

One checks that each of these exponents of $q$ is valid as a
measure of taxicab distance over its respective unit cubes.

For 
any unit cube $c_u$ in $B_{c \ge d}\thinspace ,$ 
it is valid to measure a taxicab distance by using
$(a-1) + (b-2)  + (c-1) + (d-1)\thinspace $:
its 
coordinate values $a,c,d$ are each at least 1;
and its 
coordinate 
$b$ is at least 2.

For 
any unit cube $c_u$ in $B_{c < d}\thinspace ,$
it is valid to measure a taxicab distance by using
$-1 + (a-1) + (b-2)  + (c-1) + (d-1)\thinspace $:
its 
coordinate values $a,c$ are each at least 1;
and its 
coordinates $b,d$ are at least 2.

For 
any unit cube $c_u$ in $A_{a < b}\thinspace ,$
it is valid to measure a taxicab distance by using
$(a-1) + (b-2)  + (c-1) + (d-1)\thinspace $:
its 
coordinate values $a,c,d$ are each at least 1;
and its 
coordinate $b$ is at least 2.
\\

And one sees that our summation over $B_{c\ge d}$ and 
our summation over $A_{a<b}$ are
adding the same common monomial.
We will prove ahead, using the isometry map $\psi_{z,w,y-1,x}\thinspace ,$
that the summation over $A_{a<b}$ and
the summation over $B_{c<d}$ are equal as sums of their respective monomials.
\\

We interpret the first sum of (18),
using 
taxicab distance on the unit cubes in $B_{c \ge d}.$

\begin{lemma}\label{Thm:3.3}
Let $c_u$ denote an arbitrary unit cube in block $B,$
and let the location label of $c_u$ be $(a,b,c,d).$
For the unit cubes $c_u$ in $B$ that satisfy $c \ge d,$
we have the equality
\begin{equation}
\sum_{k=1}^{n} \thinspace
               [k]_q\thinspace
               q^{k-1} 
               \frac{[k+1]_q\thinspace [k]_q}{1+q} \thinspace
= \thinspace
\sum_{c_u\in B;\thinspace c \ge d} q^{(a-1) + (b-2) + (c-1) + (d-1)}\thinspace .
\end{equation}
The exponent of $q$ on this right side is the taxicab distance between $c_u$
and the unit cube in the region $R_3(n)$ 
whose location label is $(1, 2, 1, 1).$
\end{lemma}

\begin{proof}
We re-express the sum on the right, using the definition of the block $B.$
In that definition,
when the union's index takes a given value $k$,
then for each of the $k^3$ unit cubes $c_u$ in $B$ 
which correspond to that $k,$ 
the $c_u$ has location label $(a,b,c,d)$ where
$b$ takes the value $k+1,$ and 
$a,c,d$ each take any value in $\{1, 2, \ldots k\}.$ 
Now we sum over the possible values for $k,$ 
with the further constraint that $c \ge d,$
the corresponding values of monomial
$q^{(a-1) + (b-2) + (c-1) + (d-1)} .$ 
The sum on the right side of (19) equals
\begin{equation*}
\sum_{k=1}^{n}
\Big( \sum_{1\le a \le k} q^{(a-1)} \Big)\medspace
\Big( q^{(k+1)-2} \Big)\medspace
\Big( \sum_{1 \le d \le c \le k} q^{(c-1) + (d-1)} \Big).
\end{equation*}
Within this summation:
the sum in the first pair of parentheses equals $[k]_q\thinspace .$
And the sum in the third pair of parentheses equals 
$ ( [k+1]_q\thinspace [k]_q )/ (1+q)$; 
we will prove this statement separately in the next lemma.
This completes proof of Lemma~\ref{Thm:3.3}.
\end{proof}

\begin{lemma}\label{Thm:3.4}
Let $k$ be a fixed positive integer. 
We have the equality
\begin{equation*}
\sum_{1 \le d \le c \le k}
     q^{(c-1) + (d-1)} = 
     \frac{[k+1]_q \thinspace [k]_q}{1+q}\thinspace .
\end{equation*}
\end{lemma}

\begin{proof}
Let $M$ be the matrix 
$M = [m_{i,j}]_{1 \le i \le k;\thinspace 1 \le j \le k+1}\thinspace ,$
where $m_{i,j} = q^{(i-1) + (j-1)}.$

Let $S = \sum_{1 \le j \le i \le k} q^{(i-1) + (j-1)} ,$
the sum of the entries in $M$ whose locations satisfy $i \ge j.$
We interpret $S$ as the left side of the lemma's equation, 
after a change of variables.
And let $T = \sum_{1 \le i < j \le k+1} q^{(i-1) + (j-1)} ,$
the sum of the entries in $M$ whose locations satisfy $i < j.$ 

Let $f$ be the mapping 
defined by the rule $f : m_{i,j} \mapsto m_{j, i+1}.$
One checks that $f$ is a one-to-one correspondence from 
the set of all entries $m_{i,j}$ in $M$ 
which satisfy $1 \le j \le i \le k,$
to the set of all entries $m_{i,j}$ in $M$ 
which satisfy $1 \le i < j \le k+1.$
And one checks that for each $m_{i,j}$ 
in this set of pre-image elements for $f,$
the image value satisfies that
$f(m_{i,j}) = q^1 \cdot m_{i,j}.$
So there follows by summing: $\thinspace T = f(S) = q^1 \cdot S.$ 

Let $U$ be the sum total of all of the entries in $M,$ so
$U = \sum_{1 \le i \le k;\thinspace
          1 \le j \le k+1} q^{i-1}\thinspace q^{j-1} .
$
We see by algebra, on the sum of these entries in this rectangular matrix:
\begin{equation*}
  U\hspace*{0.250em}
    =\hspace*{0.250em} 
      \sum_{1 \le i \le k;\thinspace
          1 \le j \le k+1} q^{i-1} \cdot q^{j-1}\hspace*{0.250em}
    =\hspace*{0.250em} 
       \big(\sum_{1 \le i \le k  } q^{i-1}\big)
       \cdot
       \big(\sum_{1 \le j \le k+1} q^{j-1}\big)\hspace*{0.250em}
    =\hspace*{0.250em} 
    [k]_q\thinspace
       [k+1]_q .   
\end{equation*}
And by our partition of the set of all locations in matrix $M,$ 
we have $U = S + T .$
So by substitution: $S + q^1 \cdot S = [k]_q\thinspace [k+1]_q .$
Conclude by dividing each side by $(1+q).$
\end{proof}

\par\noindent\rule{\textwidth}{0.4pt}

\vspace*{2.50ex}

See the Figure 4 for illustration. 
\\



\begin{figure}[h] 
    \centering

    \includegraphics[scale=1.0]{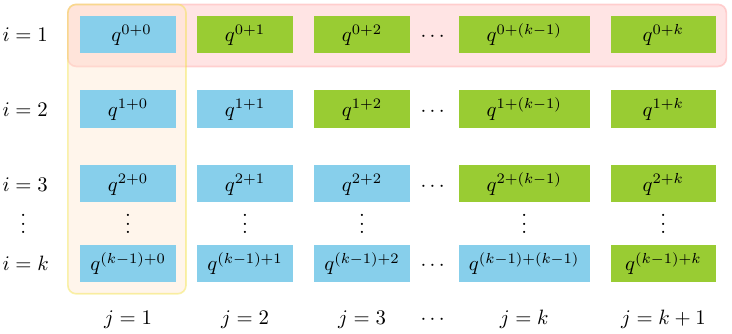}
    \caption{
    To prove the equation
    $\sum_{1 \le d \le c \le k} q^{(c-1) + (d-1)} 
     =\thinspace ( [k+1]_q\thinspace [k]_q)/(1+q) .$
    }
\end{figure}



Now we interpret the second summation of (18).
Our plan is to use the remaining part of our partition of block $B;$
and then map that onto 
the 
subcomponent of block $A.$

\begin{lemma}\label{Thm:3.5}
Let $c_u$ denote an arbitrary unit cube in block $B,$
and let the location label of $c_u$ be $(a,b,c,d).$
For the unit cubes $c_u$ in $B$ that satisfy $c < d,$
we have the equality
\begin{equation}
\sum_{k=1}^{n} \thinspace
               [k]_q\thinspace
               q^{k-1} 
               \frac{[k]_q\thinspace [k-1]_q}{1+q} \thinspace
= \thinspace
q^{-1} \cdot  \sum_{c_u\in B;\thinspace c < d} 
  q^{(a-1) + (b-2) + (c-1) + (d-1)}\thinspace .
\end{equation}
\end{lemma}

\begin{proof}
Start with the same reasoning as used in 
the start of our proof of Lemma~\ref{Thm:3.3} 
to re-express the sum on the right side of (19),
but now applied to the unit cubes $c_u$ of $B$
which 
satisfy $c < d .$ We obtain that
\begin{equation*}
\sum_{k=1}^{n}
\Big( \sum_{1\le a \le k} q^{(a-1)} \Big)\thinspace
\Big( q^{(k+1)-2} \Big)\thinspace
\Big( \sum_{1 \le c < d \le k} q^{(c-1) + (d-1)} \Big)\thinspace
=\thinspace
\sum_{c_u\in B;\thinspace c < d} q^{(a-1) + (b-2) + (c-1) + (d-1)} .
\end{equation*}

One checks that the sum
$\sum_{1 \le c < d \le k} q^{(c-1) + (d-1)}$ 
equals
$q^{1} \cdot (\thinspace [k]_q\thinspace [k-1]_q\thinspace )\thinspace 
 / \thinspace (1+q)\thinspace .$
That is: by applying the change of variable $d^{\prime} = d-1 ,$
and by Lemma~\ref{Thm:3.4}, 
we have
\begin{equation*}
\sum_{1 \le c < d \le k} q^{(c-1) + (d-1)} \thinspace 
 = \medspace
q^{1} \cdot \sum_{1 \le c \le d^{\prime} \le k-1} 
   q^{(c-1) + (d^{\prime}-1)} \thinspace 
 = \thinspace 
q^{1} \cdot \frac{[k]_q\thinspace [k-1]_q}{1+q}\thinspace .
\end{equation*}
Now substitute this into the equation in the first paragraph, 
and multiply each side by $q^{-1},$ 
to complete the proof of Lemma~\ref{Thm:3.5}. 
\end{proof}

Now for this subcomponent $B_{c < d}\thinspace ,$
we apply the isometry mapping $\psi_{z,w,y-1,x}$ from (16).
This $\psi_{z,w,y-1,x}$ maps the point $(1,2,1,1)$ to $(1,1,1,1).$
So as an isometry in $\mathbb{R}^4,$
we have that for any unit cube $c_u$ in block $B,$
the taxicab distance between $c_u$ and $(1,2,1,1)$ is equal to 
the taxicab distance between $\psi_{z,w,y-1,x}(c_u)$ and $(1,1,1,1).$
We obtain a new equation for the right side of (20).
\begin{equation}
q^{-1}\cdot  \sum_{c_u \in B;\thinspace c<d} q^{(a-1) + (b-2) + (c-1) + (d-1)}
=\thinspace
q^{-1}\cdot  \sum_{c_u \in A;\thinspace a  < b } 
      q^{(a -1) + (b -1) + (c -1) + (d -1)}
\end{equation}
 
We apply this to Lemma~\ref{Thm:3.5}, 
to obtain a new interpretation of the second summation of (18).
Here we consolidate the product 
$q^{-1} \cdot q^{(b - 1)} .$

\begin{lemma}\label{Thm:3.6}
Let $c_u$ denote an arbitrary unit cube in block $A,$
and let the location label of $c_u$ be $(a,b,c,d).$
For the unit cubes $c_u$ of $A$
that satisfy $a<b,$ we have the equality
\begin{equation}
\sum_{k=1}^{n}\thinspace
               [k]_q\thinspace
               q^{k-1}\thinspace 
               \frac{[k]_q\thinspace [k-1]_q}{1+q} \thinspace
= \thinspace
\sum_{c_u\in A;\thinspace a < b} q^{(a-1) + (b-2) + (c-1) + (d-1)}\thinspace .
\end{equation}
The exponent of $q$ on this right side is the taxicab distance 
between $c_u$ and the unit cube in the region $R_3(n)$ 
whose location label is $(1,2,1,1).$\\
\end{lemma}

Lemmas~\ref{Thm:3.3} and~\ref{Thm:3.6} 
each give an interpretation of 
one of the two summations of the left side of equation (18).
We conclude by interpreting their sum,
and demonstrate that it equals the right side of (18).

\begin{lemma}\label{Thm:3.7}
Let $c_u$ denote an arbitrary unit cube sized location
in the region $R_3(n),$ and 
let the location label of $c_u$ be $(x,y,z,w) = (a,b,c,d).$
We have the equalities
\begin{equation*}
\begin{split}
\sum_{c_u\in B;\thinspace c \ge d}
      q^{(a-1) + (b-2) + (c-1) + (d-1)} \thinspace
 + \thinspace
\sum_{c_u\in A; \thinspace a < b}
      q^{(a-1) + (b-2) + (c-1) + (d-1)} \thinspace
&= \\
\sum_{\substack{
                c_u\in R_3(n)\\
                1\le a < b\le n+1 \\
                1 \le d \le c\le n 
               }
      }
      q^{(a-1) + (b-2) + (c-1) + (d-1)} \thinspace 
&= \thinspace  
{\genfrac{[}{]}{0pt}{0}{n+1}{2}}_q^2
\thinspace .
\end{split}
\end{equation*}
\end{lemma}

\begin{proof}
The first equation follows from Lemma 14, for
the union of the disjoint sets 
$B_{c\ge d}$ and $A_{a<b} \thinspace .$ 
To establish the second equation: 
we re-express the summation by using algebra, 
and by the definition of the rectangular region $R_3(n) ,$
to obtain
\begin{equation*} 
\Big(\sum_{1 \le a < b \le n+1} q^{(a-1) + (b-2)}\Big) \thinspace
\Big(\sum_{1 \le d \le c \le n} q^{(c-1) + (d-1)}\Big) \thinspace .
\end{equation*}
Each of these final two sums equals 
$
{\genfrac{[}{]}{0pt}{1}{n+1}{2}}_q
,$
by Lemma~\ref{Thm:3.4} 
and by applying the change of variable $b^{\prime} = b-1 .$ 
Our equality is now established.
\end{proof}

\noindent
Summary:

We have constructed 3 summations: 
over the unit cubes $c_u$ in $B_{c\ge d}\thinspace ,$ and
in $B_{c < d}\thinspace ,$ and
in $A_{a < b}\thinspace .$
Each of these summations is of some monomial 
that is defined in terms of 
a $q$-counting of the taxicab distance
between $c_u$ and the start point $(1,2,1,1).$

Our summations over $B_{c\ge d}\thinspace$ and $A_{a < b}\thinspace $
are both of a common monomial. 
Since the union of the disjoint sets 
$B_{c\ge d}\thinspace$ and $A_{a < b}\thinspace $
is the solid $S_{3,3}(n) ,$
we obtain that their summation total
equals the $q$-analogue of ${\binom{n+1}{2}}^2$ 
on the right side of (18).

Our summations over $B_{c < d}\thinspace$ and $A_{a < b}\thinspace $
are equal, 
due to isometry $\psi_{z, w, y-1, x}\thinspace .$
Since the union of the disjoint sets 
$B_{c\ge d}\thinspace$ and $B_{c < d}\thinspace $
is the block $B ,$ 
we obtain that their summation total
equals the $q$-analogue of $\sum_{i=1}^{n} i^3$ 
on the left side of (18).
\\ 

This completes our interpretation of Garrett and Hummel's $q$-analogue equation. 
\\


\subsection{Interpreting the q-analogue from Warnaar}
In this section we offer \emph{two} interpretations of
the $q$-analogue of the sum of third powers
equation from Warnaar~\cite{example11_warnaar}. 

The equation of Warnaar states:
\begin{equation}
\sum_{k=1}^{n}\thinspace
  q^{2n-2k}\thinspace
  \frac{ (1-q^k)^2 (1-q^{2k})} {(1-q)^2 (1-q^2)}\thinspace
=\thinspace
{ { \genfrac{[}{]}{0pt}{0}{n+1}{2} }
}_q^2\thinspace .
\end{equation}
To enable interpretation with taxicab distance, 
we re-express the left side with algebra:
first as
\begin{equation*}
\sum_{k=1}^{n}\thinspace 
(q^{n-k})^2\thinspace 
([k]_q)^2 \thinspace 
\frac{(1 - q^{k+1}) +\thinspace  q^{k+1} (1- q^{k-1})}{1-q^2}\thinspace ,
\end{equation*}
and then in the form
\begin{equation}
\begin{split}
&\sum_{k=1}^n\thinspace 
  (q^{n-k}\thinspace [k]_q) \cdot 
  (q^{n-k}) \cdot 
  \frac{[k+1]_q\thinspace [k]_q}{1+q}\\
+ 
&\sum_{k=1}^n\thinspace 
  [k]_q \cdot 
  q^{k-1} \cdot 
  q^{2} 
        \frac{ (q^{n-k}\thinspace[k]_q)\thinspace 
               (q^{n-k}\thinspace[k-1]_q)}{1+q} 
\thinspace
=
\thinspace
{ { \genfrac{[}{]}{0pt}{0}{n+1}{2} }
}_q^2\thinspace .
\end{split}
\end{equation}

We begin our first interpretation of (24) by selecting the block $B .$
Let $c_u$ denote an arbitrary unit cube in $B ,$ and let its location label 
be $(x,y,z,w) = (a,b,c,d) .$
Partition $B$ into the parts $B_{c \ge d}$ and $B_{c < d}\thinspace .$ 

For this $q$-analogue equation, we will select for our $q$-counts of taxicab distances:
\begin{itemize}
\item
$\hspace*{0.250em}\hphantom{q^{1}\cdot }\sum_{c_u \in B;\thinspace c \ge d}\thinspace
 q^{(n-a) + (n+1-b) + (c-1) + (d-1)}$

\item
$q^{1}\cdot \sum_{c_u \in B;\thinspace c < d}\thinspace
 q^{(a-1) + (b-2) + (n-c) + (n-d)}$

\item
$q^{1}\cdot \sum_{c_u \in A;\thinspace a < b}\thinspace
 q^{(n-a) + (n-b) + (c-1) + (d-1)}\thinspace
 =\thinspace
 \sum_{c_u \in A;\thinspace a < b}\thinspace
 q^{(n-a) + (n+1-b) + (c-1) + (d-1)}\thinspace .
 $
\end{itemize}

One checks that each of these exponents of $q$ is valid as a
measure of taxicab distance over its respective unit cubes.

For 
any unit cube $c_u$ in $B_{c \ge d}\thinspace ,$ 
it is valid to measure a taxicab distance by using
$(n-a) + (n+1-b)  + (c-1) + (d-1)\thinspace $:
its 
coordinate value $a$ is at most $n ;$
and its 
coordinate $b$ is at most $n+1 ;$ 
and its 
coordinates $c,d$ are each at least 1.

For 
any unit cube $c_u$ in $B_{c < d}\thinspace ,$
it is valid to measure a taxicab distance by using
$1 + (a-1) + (b-2)  + (n-c) + (n-d)\thinspace $:
its 
coordinate value $a$ is at least 1;
and its 
coordinate $b$ is at least 2; 
and its 
coordinate $c$ is at most $n-1$;
and its 
coordinate $d$ is at most $n .$

For 
any unit cube $c_u$ in $A_{a < b}\thinspace ,$
it is valid to measure a taxicab distance by using
$1 + (n-a) + (n-b)  + (c-1) + (d-1)\thinspace $:
its 
coordinate value $a$ is at most $n-1 ;$
and its 
coordinate $b$ is at most $n ;$ 
and its 
coordinates $c,d$ are each at least 1.
\\

And one sees that our summation over $B_{c\ge d}$ and
our summation over $A_{a < b}$ are adding the same common monomial.
We will prove ahead, using the isometry map $\psi_{z, w, y-1, x}\thinspace ,$
that the summation over $A_{a< b}$ and 
the summation over $B_{c < d}$ are equal as sums of their respective monomials.
\\
 
We interpret the first sum of (24), using taxicab distance on the unit
cubes in $B_{c \ge d} .$

\begin{lemma}\label{Thm:3.8}
Let $c_u$ denote an arbitrary unit cube in block $B ,$
and let the location label of $c_u$ be $(a,b,c,d) .$
For the unit cubes $c_u$ in $B$ that satisfy $c\ge d ,$
we have the equality
\begin{equation}
\sum_{k=1}^n\thinspace 
  (q^{n-k}\thinspace [k]_q) \cdot 
  (q^{n-k}) \cdot \frac{[k+1]_q\thinspace [k]_q}{1+q}\thinspace 
= \thinspace
\sum_{c_u \in B;\thinspace c \ge d}
 q^{(n-a) + (n+1-b) + (c-1) + (d-1)}\thinspace .
\end{equation}
The exponent of $q$ on this right side is the taxicab distance between
$c_u$ and the unit cube in the region $R_3(n)$ whose location label is
$(n, n+1, 1, 1) .$
\end{lemma}
\begin{proof}
We re-express the sum 
on the right, using the definition of block $B ;$
we follow the same reasoning as in the proof of Lemma~\ref{Thm:3.3}. 
The sum on the right side of (25) equals
\begin{equation*}
\sum_{k=1}^{n}\thinspace 
    \Big( \sum_{1 \le a \le k} q^{(n-a)} \Big) \thinspace 
    \Big( q^{n+1 - (k+1)} \Big) \thinspace 
    \Big( \sum_{1 \le d \le c \le k} q^{(c-1)+(d-1)} \Big)\thinspace .
\end{equation*}
The sum in the first pair of parentheses equals 
$q^{n-k}\thinspace [k]_q \thinspace .$
For the sum in the third pair of parentheses, 
we apply Lemma~\ref{Thm:3.4}. 
This completes proof of Lemma~\ref{Thm:3.8}. 
\end{proof}

Now we interpret the second summation of (24).
\begin{lemma}\label{Thm:3.9}
Let $c_u$ denote an arbitrary unit cube in block $B ,$
and let the location label of $c_u$ be $(a,b,c,d) .$
For the unit cubes $c_u$ in $B$ that satisfy $c < d ,$
we have the equality
\begin{equation}
\sum_{k=1}^n\thinspace 
  [k]_q \cdot 
  q^{k-1} \cdot 
  q^{2} 
        \frac{ (q^{n-k}\thinspace[k]_q)\thinspace 
               (q^{n-k}\thinspace[k-1]_q)}{1+q} 
\thinspace = \thinspace
q^{1} \cdot \sum_{c_u \in B;\thinspace c < d}
 q^{(a-1) + (b-2) + (n-c) + (n-d)}\thinspace .
\end{equation}
\end{lemma}
\begin{proof}
We re-express the sum on the right side of (26), 
using the definition of block $B .$
We obtain that
\begin{equation*}
\sum_{k=1}^{n}\thinspace 
    \Big( \sum_{1 \le a \le k} q^{(a-1)} \Big) \thinspace 
    \Big( q^{(k+1)-2} \Big) \thinspace 
    \Big( 
          \sum_{1 \le c < d \le k} q^{(n-c)+(n-d)} \Big) 
    \thinspace
=\thinspace\sum_{c_u \in B;\thinspace c < d} 
 q^{(a-1) + (b-2) + (n-c) + (n-d)}
\thinspace .
\end{equation*}
And now for the sum in the third pair of parentheses, we apply the
change of variables 
$c^{\prime} = k+1-c ,$ and 
$d^{\prime} = k+1-d ,$ and our proof of Lemma~\ref{Thm:3.5}. 
We obtain that
\begin{equation*}
\sum_{1 \le c < d \le k} q^{(n-c)+(n-d)} 
=\thinspace
q^{n-k} q^{n-k} 
  \sum_{k \ge c^{\prime} > d^{\prime} \ge 1} 
  q^{(c^{\prime}-1)+(d^{\prime}-1)} 
=\thinspace
q^{n-k} q^{n-k} \cdot q^{1}\cdot 
  \frac{[k]_q\thinspace [k-1]_q}{1+q}\thinspace . 
\end{equation*}
Now substitute this into the 
first displayed equation in this proof, and
multiply each side by $q^1 ,$ to complete the proof of Lemma~\ref{Thm:3.9}. 
\end{proof}

Now for this subcomponent $B_{c<d}\thinspace ,$ 
we apply the isometry mapping $\psi_{z, w, y-1, x}$ from (16).
This $\psi_{z, w, y-1, x}$ maps the point $(1,2,n,n)$ to 
$(n,n,1,1) .$
So for any unit cube $c_u$ in block $B ,$ 
the taxicab distance between $c_u$ and
$(1,2,n,n)$ is equal to the taxicab distance between
$\psi_{z, w, y-1, x}(c_u)$ and $(n,n,1,1).$ 
We obtain a new equation for the right side of (26).
\begin{equation}
q^{1} \cdot  \sum_{c_u \in B;\thinspace c < d}
 q^{(a-1) + (b-2) + (n-c) + (n-d)}\thinspace
=\medspace
q^{1} \cdot \sum_{c_u \in A;\thinspace a < b}
 q^{(n-a) + (n-b) + (c-1) + (d-1)}\thinspace .
\end{equation}

We apply this to Lemma~\ref{Thm:3.9}, 
to obtain a new interpretation of the second
summation of (24). Here we consolidate the product $q^1\cdot q^{(n-b)} .$

\begin{lemma}\label{Thm:3.10}
Let $c_u$ denote an arbitrary unit cube in block $A ,$
and let the location label of $c_u$ be $(a,b,c,d) .$
For the unit cubes $c_u$ in $A$ that satisfy $a < b ,$
we have the equality
\begin{equation}
\sum_{k=1}^n\thinspace 
  [k]_q \cdot 
  q^{k-1} \cdot 
  q^{2} 
        \frac{ (q^{n-k}\thinspace[k]_q)\thinspace 
               (q^{n-k}\thinspace[k-1]_q)}{1+q} 
\thinspace = \thinspace
\sum_{c_u \in A;\thinspace a < b}
 q^{(n-a) + (n+1-b) + (c-1) + (d-1)}\thinspace .
\end{equation}
The exponent of $q$ on this right side is the taxicab distance between
$c_u$ and the unit cube in the region $R_3(n)$ whose location label
is $(n, n+1, 1, 1) .$
\end{lemma}

Lemmas~\ref{Thm:3.8} and~\ref{Thm:3.10} 
each give an interpretation of 
one of the two summations of the left side of equation (24).
We conclude by interpreting their sum, 
and demonstrate that it equals the right side of (24).

\begin{lemma}\label{Thm:3.11}
Let $c_u$ denote an arbitrary unit cube sized location in the region $R_3(n) ,$
and let the location label of $c_u$ be $(x,y,z,w) = (a,b,c,d) .$
We have the equalities
\begin{equation*}
\begin{split}
\sum_{c_u \in B;\thinspace c \ge d}
 q^{(n-a) + (n+1-b) + (c-1) + (d-1)}\thinspace
+\thinspace\sum_{c_u \in A;\thinspace a < b}
 q^{(n-a) + (n+1-b) + (c-1) + (d-1)}\thinspace
 &=\\
\sum_{\substack{c_u \in R_3(n);\\
                1\le a < b \le n+1\\
                1\le d \le c \le n
               }
     }
 q^{(n-a) + (n+1-b) + (c-1) + (d-1)}\thinspace
 &=\thinspace {\genfrac{[}{]}{0pt}{0}{n+1}{2}}_q^2\thinspace .
\end{split}
\end{equation*}
\end{lemma}
\begin{proof}
The first equation follows from Lemma 14, for
the union of the disjoint sets 
$B_{c\ge d}$ and $A_{a<b}\thinspace .$
To establish the second equation: 
we re-express the summation by using algebra, 
and by the definition of the rectangular region $R_3(n) ,$ to obtain
\begin{equation*} 
 \Big(\sum_{1\le a < b \le n+1} q^{(n-a) + (n+1-b)}\Big)
  \thinspace \Big(\sum_{1\le d \le c \le n} q^{(c-1) + (d-1)}\Big)\thinspace .
\end{equation*}
Each of these final two sums equals 
${\genfrac{[}{]}{0pt}{1}{n+1}{2}}_q\thinspace .$
For the sum in the first pair of parentheses,
apply the change of variables
$\thinspace a^{\prime} = n+2-a\thinspace $ and 
$\thinspace b^{\prime} = n+2-b,$
to obtain that
\begin{equation*}
\sum_{1\le a < b \le n+1} q^{(n-a) + (n+1-b)}\thinspace = \thinspace
\sum_{n+1\ge a^{\prime} > b^{\prime} \ge 1} 
    q^{(a^{\prime}-2) + (b^{\prime}-1)}
\thinspace .
\end{equation*}
This has the same form as in the proof of Lemma~\ref{Thm:3.7}.
Our equality is established.
\end{proof}

This completes our first interpretation of Warnaar's $q$-analogue equation. 
\\

For our second interpretation of Warnaar's $q$-analogue equation, 
we conduct the same steps as in our first interpretation of (24),
but now with the following changes.
\begin{itemize}
\item 
We begin with the original equation from Warnaar as stated in (23).
But now we re-express the left side with algebra, first as
\begin{equation*}
\begin{split}
&\sum_{k=1}^{n}\thinspace 
  (q^{n-k})^2\thinspace 
  ([k]_q)^2\thinspace
  \frac{ q^{k-1} (1-q^{k+1})\thinspace +\thinspace (1 - q^{k-1}) }{1-q^2}\thinspace ,\\
&\hspace*{8.0em}\textrm{ not }\thinspace
\sum_{k=1}^{n}\thinspace 
  (q^{n-k})^2\thinspace 
  ([k]_q)^2\thinspace
  \frac{ (1-q^{k+1})\thinspace +\thinspace q^{k+1} (1 - q^{k-1}) }{1-q^2}\thinspace .
\end{split}  
\end{equation*}
And then in the form
\begin{equation}
\begin{split}
&\sum_{k=1}^{n}\thinspace 
  [k]_q\cdot 
  q^{k-1}\cdot 
  \frac{ (q^{n-k}\thinspace [k+1]_q)\thinspace (q^{n-k}\thinspace[k]_q) }{1+q}\\
&\hspace*{5.00em}+\thinspace \sum_{k=1}^{n}\thinspace 
  q^{n-k}\thinspace [k]_q\cdot
  q^{n-k}\cdot
  \frac{ [k]_q \thinspace [k-1]_q }{1+q}\thinspace 
=\thinspace {\genfrac{[}{]}{0pt}{0}{n+1}{2}}_q^{2}\thinspace ,
\end{split}  
\end{equation}
not
\begin{equation*}
\begin{split}
&\sum_{k=1}^{n}\thinspace 
  (q^{n-k}\thinspace [k]_q)\cdot 
  (q^{n-k})\cdot 
  \frac{ [k+1]_q \thinspace [k]_q }{1+q}\\
&\hspace*{5.00em}+\thinspace \sum_{k=1}^{n}\thinspace 
  [k]_q\cdot
  q^{k-1}\cdot
  q^2 \frac{(q^{n-k}\thinspace [k]_q) \thinspace 
            (q^{n-k}\thinspace [k-1]_q) }{1+q}\thinspace 
=\thinspace {\genfrac{[}{]}{0pt}{0}{n+1}{2}}_q^{2}\thinspace .
\end{split}  
\end{equation*}
\item 

For this $q$-analogue equation, we will 
select for our $q$-counts of taxicab distances:

\begin{itemize}
\item[\ding{70}]
$\thinspace \hspace*{0.250em}\hphantom{q^{-1}\cdot }\sum_{c_u \in B;\thinspace c \ge d}\thinspace
 q^{(a-1) + (b-2) + (n-c) + (n-d)}$

\item[\ding{70}]
$\thinspace q^{-1}\cdot \sum_{c_u \in B;\thinspace c < d}\thinspace
 q^{(n-a) + (n+1-b) + (c-1) + (d-1)}$

\item[\ding{70}]
$\thinspace q^{-1}\cdot \sum_{c_u \in A;\thinspace a < b}\thinspace
 q^{(a-1) + (b-1) + (n-c) + (n-d)}\thinspace
 =\thinspace
 \sum_{c_u \in A;\thinspace a < b}\thinspace
 q^{(a-1) + (b-2) + (n-c) + (n-d)}\thinspace .
 $
\end{itemize}

One checks that each of these exponents of $q$ is valid as a
measure of taxicab distance over its respective unit cubes.

For 
any unit cube $c_u$ in $B_{c \ge d}\thinspace ,$ 
it is valid to measure a taxicab distance by using
$(a-1) + (b-2)  + (n-c) + (n-d)\thinspace $:
its 
coordinate value $a$ is at least 1;
and its 
coordinate $b$ is at least 2; 
and its 
coordinates $c,d$ are each at most $n .$

For 
any unit cube $c_u$ in $B_{c < d}\thinspace ,$
it is valid to measure a taxicab distance by using
$-1 + (n-a) + (n+1-b)  + (c-1) + (d-1)\thinspace $:
its 
coordinate value $a$ is at most $n ;$
and its 
coordinate $b$ is at most $n+1 ;$ 
and its 
coordinate $c$ is at least 1;
and its 
coordinate $d$ is at least 2.

For 
any unit cubes $c_u$ in $A_{a < b}\thinspace ,$
it is valid to measure a taxicab distance by using
$(a-1) + (b-2)  + (n-c) + (n-d)\thinspace $:
its 
coordinate value $a$ is at least 1;
and its 
coordinate $b$ is at least 2;
and its 
coordinates $c,d$ are each at most $n .$

\item 
Our interpretation of this first sum of (29) uses
\begin{equation}
\sum_{c_u \in B; c\ge d} q^{(a-1) + (b-2) + (n-c) + (n-d)}\thinspace ,
\textrm{ not }
\sum_{c_u \in B; c\ge d} q^{(n-a) + (n+1-b) + (c-1) + (d-1)}\thinspace .
\end{equation}
Here in the proof of the equality, 
we re-express this sum using the definition of block $B ,$ 
to obtain that it equals
\begin{equation*}
\sum_{k=1}^{n}\thinspace 
    \Big( \sum_{1 \le a \le k} q^{(a-1)} \Big) \thinspace 
    \Big( q^{(k+1)-2} \Big) \thinspace 
    \Big( 
          \sum_{1 \le d \le c \le k} q^{(n-c)+(n-d)} \Big) 
\thinspace .
\end{equation*}
And for this sum in the third pair of parentheses, we apply the
change of variables 
$c^{\prime} = k+1-c ,$ and 
$d^{\prime} = k+1-d ,$ and 
Lemma~\ref{Thm:3.4}. 
We obtain that
\begin{equation*}
\sum_{1 \le d \le c \le k} q^{(n-c)+(n-d)} 
=\thinspace
q^{n-k} q^{n-k} 
  \sum_{k \ge d^{\prime} \ge c^{\prime} \ge 1} 
  q^{(c^{\prime}-1)+(d^{\prime}-1)} 
=\thinspace
q^{n-k} q^{n-k}\thinspace 
  \frac{[k+1]_q\thinspace [k]_q}{1+q}\thinspace , 
\end{equation*}
to complete the proof that this summation equals the
first sum of (29).

\item
Our initial interpretation of this second sum of (29) uses
\begin{equation}
q^{-1}\cdot \sum_{c_u \in B; c < d} q^{(n-a) + (n+1-b) + (c-1) + (d-1)}\thinspace ,
\textrm{ not }
\thinspace q^{1}\cdot \sum_{c_u \in B; c < d} q^{(a-1) + (b-2) + (n-c) + (n-d)}\thinspace .
\end{equation}
\item
Our interpretation of this second sum of (29),  
after we apply the isometry mapping
$\psi_{z, w, y-1, x}\thinspace ,$ uses
\begin{equation}
\sum_{c_u \in A; a < b} q^{(a-1) + (b-2) + (n-c) + (n-d)}\thinspace ,
\textrm{ not }
\sum_{c_u \in A; a < b} q^{(n-a) + (n+1-b) + (c-1) + (d-1)}\thinspace .
\end{equation}
This $\psi_{z, w, y-1, x}$ maps the point
$(n,n+1,1,1)$ to $(1, 1, n, n) .$
So for any unit cube $c_u$ in block $B ,$
the taxicab distance between $c_u$ and $(n ,n+1, 1, 1)$
is equal to the taxicab distance between
$\psi_{z, w, y-1, x}(c_u)$ and $(1, 1, n, n) .$
\item 
Our interpretation of the total, after adding the summation
over $c_u\in B_{c\ge d}$ plus the summation over
$c_u \in A_{a<b}\thinspace ,$ uses
\begin{equation}
\sum_{\substack{c_u \in R_3(n)\\
                1\le a < b\le n+1\\
                1\le d \le c \le n}
     } 
q^{(a-1) + (b-2) + (n-c) + (n-d)}\thinspace , \textrm{ not }
\sum_{\substack{c_u \in R_3(n)\\
                1\le a < b\le n+1\\
                1\le d \le c \le n}
     }
\thinspace q^{(n-a) + (n+1-b) + (c-1) + (d-1)}\thinspace .
\end{equation}
Here to complete the proof,
we re-express the summation in (33) by using algebra,
and by the definition of the rectangular region $R_3(n) ,$ to obtain
\begin{equation*}
 \Big(\sum_{1\le a < b \le n+1} q^{(a-1) + (b-2)}\Big)
  \thinspace \Big(\sum_{1\le d \le c \le n} q^{(n-c) + (n-d)}\Big)\thinspace .
\end{equation*}
Each of these 
two sums equals 
${\genfrac{[}{]}{0pt}{1}{n+1}{2}}_q\thinspace ,$
by our proof of Lemma~\ref{Thm:3.7} and 
our proof of equality for the first sum of (29).
\end{itemize}

This completes our second interpretation of the $q$-analogue equation from 
Warnaar. 
\\

These two interpretations of Warnaar's $q$-analogue equation 
are distinct, but they have complementary 
properties with respect their use of the taxicab distance start points
$(n, n+1, 1, 1)$ and $(1, 2, n, n) .$
\\


\subsection{Interpreting the q-analogue from Zhao and Feng}
In this section we interpret the $q$-analogue of the sum of third powers
equation from Zhao and Feng~\cite{example14_zhao_etal}. 
The equation of Zhao and Feng states:
\begin{equation}
\sum_{k=1}^{n} q^{4(n-k)} 
               \frac{(1 - q^k)^2           }{(1-q  )^2}
               \frac{(1 + q^2 - 2 q^{k+1}) }{(1-q^2)  }\thinspace
=\thinspace
{ { \genfrac{[}{]}{0pt}{0}{n+1}{2} }
}_q^2\thinspace .
\end{equation}
To enable interpretation with taxicab distance, 
we re-express the left side with algebra.
\begin{equation}
\begin{split}
&\sum_{k=1}^{n}\thinspace 
(q^{n-k}\thinspace [k]_q)\cdot
(q^{n-k})\cdot
\frac{ ( q^{n-k}\thinspace [k+1]_q)\thinspace  
       ( q^{n-k}\thinspace [k  ]_q) }{1+q}\\
+\thinspace
&\sum_{k=1}^{n}\thinspace 
(q^{n-k}\thinspace [k]_q)\cdot
(q^{n-k})\cdot
q^2 \frac{ ( q^{n-k}\thinspace [k  ]_q)\thinspace 
           ( q^{n-k}\thinspace [k-1]_q) }{1+q}\thinspace
=\thinspace
{ { \genfrac{[}{]}{0pt}{0}{n+1}{2} }
}_q^2\thinspace .
\end{split}
\end{equation}

We begin our interpretation of (35) by selecting the block $B$. 
Let $c_u$ denote an arbitrary unit cube in $B ,$
and let its location label be $(x,y,z,w) = (a,b,c,d).$
Partition $B$ into the parts $B_{c \ge d}$ and $B_{c < d}\thinspace .$ 

For this $q$-analogue equation, we will select for our $q$-counts of 
taxicab distances:
\begin{itemize}
\item 
$\hspace*{0.250em}\hphantom{q^{1} \cdot }\sum_{c_u \in B;\thinspace c \ge d} \thinspace
      q^{(n-a) + (n+1-b) + (n-c) + (n-d) }$
\item 
$q^{1} \cdot  \sum_{c_u \in B;\thinspace c < d} \thinspace
       q^{(n-a) + (n+1-b) + (n-c) + (n-d) }$
\item 
$q^{1} \cdot  \sum_{c_u \in A;\thinspace a < b} \thinspace
       q^{(n-a) + (n-b) + (n-c) + (n-d) }\thinspace
=\thinspace
 \sum_{c_u \in A;\thinspace a \ge b} \thinspace
       q^{(n-a) + (n+1-b) + (n-c) + (n-d) }\thinspace .
$
\end{itemize}

One checks that each of these exponents of $q$ is valid as a
measure of taxicab distance over its respective unit cubes.

For 
any unit cube $c_u$ in $B_{c \ge d}\thinspace ,$ 
it is valid to measure a taxicab distance by using
$(n-a) + (n+1-b)  + (n-c) + (n-d)\thinspace $:
its 
coordinate values $a,c,d$ are each at most $n ;$
and its 
coordinate $b$ is at most $n+1 .$

For 
any unit cube $c_u$ in $B_{c < d}\thinspace ,$ 
it is valid to measure a taxicab distance by using
$1 + (n-a) + (n+1-b)  + (n-c) + (n-d)\thinspace $:
its 
coordinate values $a,d$ are each at most $n ;$
and its 
coordinate $b$ is at most $n+1 ;$
and its coordinate $c$ is at most $n-1$.

For 
any unit cube $c_u$ in $A_{a < b}\thinspace ,$
it is valid to measure a taxicab distance by using
$1 + (n-a) + (n-b)  + (n-c) + (n-d)\thinspace $:
its 
coordinate value $a$ is at most $n-1 ;$
and its coordinates $b,c,d$ are each at most $n .$ 
\\

And one sees that our summation over $B_{c\ge d}$ and our summation over $A_{a<b}$
are adding the same common monomial. We will prove ahead, using the isometry map
$\psi_{z, w, y-1, x}\thinspace, $ that the summation over $A_{a< b} $ 
and the
summation over $B_{c< d}$ are equal as sums of their respective monomials.
\\

We interpret the first sum of (35), using taxicab distance
on the unit cubes in $B_{c \ge d} .$

\begin{lemma}\label{Thm:3.12}
Let $c_u$ denote an arbitrary unit cube in block $B ,$
and let the location label of $c_u$ be $(a,b,c,d).$
For the unit cubes $c_u$ in $B$ that satisfy $c \ge d ,$
we have the equality
\begin{equation}
\begin{split}
\sum_{k=1}^{n}\thinspace
  (q^{n-k}\thinspace [k]_q)\thinspace \cdot\thinspace  
  (q^{n-k})\thinspace \cdot\thinspace 
  &\frac{ (q^{n-k}\thinspace [k+1]_q) \thinspace 
         (q^{n-k}\thinspace [k  ]_q) 
       }{1+q}\\
&{\hspace*{5.1250em}}=
\sum_{c_u \in B;\thinspace c \ge d}
       q^{(n-a) + (n+1-b) + (n-c) + (n-d) }\thinspace .
\end{split}
\end{equation}
The exponent of $q$ on this right side is the taxicab distance between
$c_u$ and the unit cube in the region $R_3(n)$
whose location label is $(n, n+1, n, n) .$
\end{lemma}

\begin{proof}
We re-express the sum on the right, using the definition of block $B ;$
we follow the same reasoning as in the proof of Lemma~\ref{Thm:3.3}. 
The sum on the right side of (36) equals
\begin{equation*}
\sum_{k=1}^{n}
 \Big(\sum_{1 \le a \le k} q^{(n-a)}\Big)\thinspace 
 \Big(q^{n+1 - (k+1)} \Big)\thinspace
 \Big(\sum_{1 \le d \le c \le k} q^{(n-c) + (n-d)} \Big) \thinspace .
\end{equation*}
For the sum in the third pair of parentheses, 
apply our proof of equality for the first sum of (29).
This completes proof of Lemma~\ref{Thm:3.12}. 
\end{proof}

Now we interpret the second summation of (35).
\begin{lemma}\label{Thm:3.13}
Let $c_u$ denote an arbitrary unit cube in block $B ,$
and let the location label of $c_u$ be 
$(a,b,c,d) .$
For the unit cubes $c_u$ in $B$ that satisfy $c < d ,$
we have the equality
\begin{equation}
\sum_{k=1}^n (q^{n-k}\thinspace [k]_q ) 
             \cdot
             (q^{n-k}) 
             \cdot
             q^2\thinspace 
               \frac{(q^{n-k}\thinspace [k]_q)\thinspace 
                     (q^{n-k}\thinspace [k-1]_q)}{1+q} 
\thinspace
= \thinspace                   
q^{1} \cdot \sum_{c_u \in B;\thinspace c < d}
 q^{(n-a) + (n+1-b) + (n-c) + (n-d)}\thinspace .
\end{equation}
\end{lemma}

\begin{proof}
We re-express the sum on the right side of (37),
using the definition of block $B .$
We obtain that
\begin{equation*}
\sum_{k=1}^{n}\thinspace
 \Big(\sum_{1 \le a \le k} q^{(n-a)}\Big)\thinspace 
 \Big(q^{n+1 - (k+1)} \Big)\thinspace
 \Big(\sum_{1 \le c < d \le k} 
 q^{(n-c) + (n-d)} \Big) 
\thinspace
=
\thinspace
\sum_{c_u \in B;\thinspace c < d}\thinspace
 q^{(n-a) + (n+1-b) + (n-c) + (n-d)}\thinspace .
\end{equation*}
And now 
for the sum in the third pair of parentheses,
and for the remainder of the proof,
follow the same reasoning as in our proof of Lemma~\ref{Thm:3.9}.
\end{proof}

Now for this subcomponent $B_{c < d}\thinspace ,$ 
we apply the isometry 
mapping $\psi_{z, w, y-1, x} \thinspace $ from (16).
This $\psi_{z, w, y-1, x}$ maps the point $(n, n+1, n, n)$ to $(n, n, n, n ) .$
So for any unit cube
$c_u$ in block $B ,$ the taxicab distance between $c_u$ and $(n, n+1, n, n)$
is equal to the taxicab distance between $\psi_{z, w, y-1, x}(c_u)$ and 
$(n, n, n, n) .$
We obtain a new equation for the right side of (37).
\begin{equation}
q^{1} \cdot \sum_{c_u \in B;\thinspace c < d}\thinspace
 q^{(n-a) + (n+1-b) + (n-c) + (n-d)} 
\thinspace
=\thinspace
q^{1} \cdot  \sum_{c_u \in A;\thinspace a < b}\thinspace
 q^{(n-a) + (n-b) + (n-c) + (n-d)}
\end{equation}

We apply this to Lemma~\ref{Thm:3.13}, 
to obtain new interpretation 
of the second summation of (35).
Here we consolidate the product $q^{1}\cdot q^{(n-b)} .$

\begin{lemma}\label{Thm:3.14}
Let $c_u$ denote an arbitrary unit cube in block $A ,$
and let the location label of $c_u$ be 
$(a, b, c, d) .$
For the unit cubes $c_u$ in $A$ that satisfy $a < b , $
we have the equality
\begin{equation}
\sum_{k=1}^n \thinspace
             (q^{n-k}\thinspace [k]_q ) 
             \cdot
             (q^{n-k})
             \cdot
             q^2\thinspace 
               \frac{(q^{n-k}\thinspace [k]_q)\thinspace 
                     (q^{n-k}\thinspace [k-1]_q)}{1+q} 
\thinspace
= \thinspace
\sum_{c_u \in A;\thinspace a < b}\thinspace
 q^{(n-a) + (n+1-b) + (n-c) + (n-d)}\thinspace .
\end{equation}
The exponent of $q$ on this right side is the taxicab distance between
$c_u$ and the unit cube in the region $R_3(n)$ 
whose location label is $(n, n+1, n, n) .$
\end{lemma}

Lemmas~\ref{Thm:3.12} and~\ref{Thm:3.14} 
each give an interpretation of 
one of the two summations of the left side of equation (35).
We conclude by interpreting their sum, and demonstrate that it equals
the right side of (35).

\begin{lemma}\label{Thm:3.15}
Let $c_u$ denote an arbitrary unit cube sized location in the region $R_3(n) ,$
and let the location label of $c_u$ be $(x,y,z,w) = (a,b,c,d) .$
We have the equalities
\begin{equation*}
\begin{split}
\sum_{c_u \in B;\thinspace c \ge d}
 q^{(n-a) + (n+1-b) + (n-c) + (n-d)}\thinspace
+ 
\sum_{c_u \in A;\thinspace a < b}
 q^{(n-a) + (n+1-b) + (n-c) + (n-d)}
&=\\
\sum_{\substack{
      c_u \in R_3(n)\\
      1 \le a <   b \le n+1\\ 
      1 \le d \le c \le n
                }
     }
 q^{(n-a) + (n+1-b) + (n-c) + (n-d)}
&=\thinspace
{\genfrac{[}{]}{0pt}{0}{n+1}{2}}_q^2 \thinspace .
\end{split}
\end{equation*}
\end{lemma}

\begin{proof}
The first equation follows from Lemma 14, for
the union of the disjoint sets 
$B_{c\ge d}$ and $A_{a<b}\thinspace .$
To establish the second equation: 
we re-express the summation by using algebra, 
and by the definition of
the rectangular region $R_3(n) ,$ to obtain
\begin{equation*}
\begin{split}
\Big(
\sum_{ 1 \le a < b \le n+1
      }
 q^{(n-a) + (n+1-b)}
\Big)\thinspace
\Big(
\sum_{ 1 \le d \le c \le n }
 q^{(n-c) + (n-d)}
\Big) \thinspace .
\end{split}
\end{equation*}
Each of these two sums equals
${\genfrac{[}{]}{0pt}{1}{n+1}{2}}_q .$
For the sum in the first pair of parentheses,
see the proof of Lemma~\ref{Thm:3.11}. 
And for the sum in the second, 
see our proof of equality for the first sum of (29).
\end{proof}

This completes our interpretation of Zhao and Feng's $q$-analogue equation. 
\\

\begin{remark}\label{Thm:Remark}
For our interpretation of the $q$-analogue equation of Garrett and Hummel,
we used $q$-statistics that applied taxicab distance 
from the start point $(1,2,1,1)$ in $R_{3}(n) .$

For our two interpretations of the $q$-analogue equation of Warnaar,
we used $q$-statistics that applied taxicab distance
from the start points $(n,n+1,1,1)$ and $(1,2,n,n)$ in $R_{3}(n) .$

For our interpretation of the $q$-analogue equation of Zhao and Feng,
we used $q$-statistics that applied taxicab distance 
from the start point $(n,n+1,n,n)$ in $R_{3}(n) .$
\\

These four selected start points for taxicab distance measure
each coincide with four of the corner vertices
of our assembled solid 
$S_{3,3}(n) .$
That is:

For the step triangle in the $(x,y)$-plane which comprises the
$1 + 2 + \ldots + n$ unit squares whose location labels 
each satisfy $1 \le a < b \le n+1 ,$
the particular unit squares whose coordinates
in $(x,y)$ are both as small as possible,
or are both as large as possible, are
at $(x,y) = (1,2) $ 
and at $(x,y) = (n,n+1) .$

For the step triangle in the $(z,w)$-plane which comprises the
$1 + 2 + \ldots + n$ unit squares whose location labels 
each satisfy $1 \le d \le c \le n ,$
the 
unit squares whose coordinates 
in $(z,w)$ are both as small as possible,
or are both as large as possible, are
at $(z,w) = (1,1) $ 
and at $(z,w) = (n,n) .$
\\
\end{remark}


\subsection{Interpreting the q-analogue from 
            Forster, Garrett, Jacobsen, and Wood}
In this section we interpret the $q$-analogue of the sum of third powers
equation from Forster, Garrett, Jacobsen, 
              and Wood~\cite{example16_forster_etal}. %
Their equation states:
\begin{equation}
\sum_{i=1}^{n}\thinspace  2q^{i-1}\thinspace 
               \bigg(\frac{1 - q^{i-1}}{1-q}\bigg)^2\thinspace
               \bigg(
               \Big(\frac{1 - q^{i-2}}{1-q}\Big)
               +
               \Big(\frac{1 - q^{i}}{1-q}\Big)
               \bigg)\thinspace
=\thinspace
(2q + 2q^2)\thinspace
{ { \genfrac{[}{]}{0pt}{0}{n}{2} }
}_q^2\thinspace .
\end{equation}
We re-express this left side, first with cosmetic changes.
\begin{enumerate}
\item 
This summand, when $i=1 ,$ equals 0.
We reset the initial value of the index to $i=2.$
\item 
Perform the change of variable $j = i-1 ,$ 
then rename the index $j$ to be $i .$
\item 
For consistency with previous formulas, 
reassign the index end value. 
Perform the change of variable $N = n-1 ,$
then rename the variable $N$ to be $n .$
\end{enumerate}
\begin{equation*}
\sum_{i=1}^{n}\thinspace 2q^{i}\thinspace 
               \bigg(\frac{1 - q^{i}}{1-q}\bigg)^2\thinspace
               \bigg(
               \Big(\frac{1 - q^{i-1}}{1-q}\Big)
               +
               \Big(\frac{1 - q^{i+1}}{1-q}\Big)
               \bigg)
               \thinspace
=\thinspace
(2q + 2q^2)\thinspace
{ { \genfrac{[}{]}{0pt}{0}{n+1}{2} }
}_q^2\thinspace .
\end{equation*}
To enable interpretation with taxicab distance, 
we re-express the left side with algebra: first as
\begin{equation*}
2 q \cdot 
               \Big(
               \sum_{i=1}^{n}\thinspace  
               [i]_{q}\thinspace
               q^{i-1}\thinspace
               \big(\thinspace [i+1]_q\thinspace [i  ]_q\thinspace \big)\thinspace 
               +\thinspace 
               \sum_{i=1}^{n}\thinspace 
               [i]_{q}\thinspace
               q^{i-1}\thinspace
               \big(\thinspace [i  ]_q\thinspace [i-1]_q\thinspace \big)
                \Big)\thinspace ,
\end{equation*}
and then in the form
\begin{equation}
\begin{split}
2 q^{1} \cdot 
               \Big(
               &\sum_{i=1}^{n}\thinspace
               [i]_{q}\thinspace
               q^{i-1}\thinspace
               \frac{[i+1]_q\thinspace [i  ]_q\thinspace}{1+q} 
               \thinspace +\thinspace    \sum_{i=1}^{n}\thinspace  
               [i]_{q}\thinspace
               q^{i-1}\thinspace
               \frac{[i  ]_q\thinspace [i-1]_q\thinspace}{1+q} 
                \Big)\medspace
\\                
+\thinspace 2 q^2 \cdot 
               \Big(
               &\sum_{i=1}^{n}\thinspace  
               [i]_{q}\thinspace
               q^{i-1}\thinspace
               \frac{[i+1]_q\thinspace [i  ]_q\thinspace}{1+q} 
               \thinspace +\thinspace \sum_{i=1}^{n}\thinspace
               [i]_{q}\thinspace
               q^{i-1}\thinspace
               \frac{[i  ]_q\thinspace [i-1]_q\thinspace}{1+q} 
                \Big)\medspace
=\medspace
(2q + 2q^2)\thinspace
{ { \genfrac{[}{]}{0pt}{0}{n+1}{2} }
}_q^2\thinspace .
\end{split}
\end{equation}
Each of the two parenthesized expressions on the left side of this equation
is identical to the left side of Garrett and Hummel's $q$-analogue
as expressed in the form (18),
\begin{equation*}
               \sum_{k=1}^{n}\thinspace  
               [k]_{q}\thinspace
               q^{k-1}\thinspace
               \frac{[k+1]_q\thinspace [k  ]_q}{1+q}
               +\thinspace 
               \sum_{k=1}^{n}\thinspace  
               [k]_{q}\thinspace
               q^{k-1}\thinspace
               \frac{[k  ]_q\thinspace [k-1]_q}{1+q} 
               \medspace
=\medspace
{ { \genfrac{[}{]}{0pt}{0}{n+1}{2} }
}_q^2\thinspace ,
\end{equation*}
after changing the index of summation from $k$ to $i .$
\\

Our plan for interpretation of (41) is:
\begin{enumerate}
\item
We interpret the first parenthesized expression on the left side of (41), 
incorporating its outside factor $q^1 ,$
in a manner similar to our interpretation of (18).
\item
We interpret the second parenthesized expression on left side of (41), 
incorporating its outside factor $q^2 ,$
in another similar manner 
but using sets of unit cubes which are disjoint from those
used for interpreting the first parenthesized expression.
\item
Then
we interpret the entire equation (41)
by using the union of all of these disjoint sets of unit cubes, 
to form a reassembled solid in $R_3(n) .$ 
\end{enumerate}
For 
our interpretation of (41)
we assign no significance to the outside factor of $2 ,$
other than that we would perform our counting procedures twice 
to obtain totals in that form. 
\\

Here we begin 
step 1 of our interpretation of equation (41):
for the first parenthesized expression on the left side, 
incorporating its outside factor $q^1 .$
\\

Select the block $B.$
Let $c_u$ denote an arbitrary unit cube in $B,$ and let its location label
be $(x,y,z,w) = (a,b,c,d) .$ 
Partition $B$ into 
the parts
$B_{c \ge d}$ 
and $B_{c < d}\thinspace .$ 

For this $q$-analogue equation, we will select for our $q$-counts of 
taxicab distances:
\begin{itemize}
\item 
$q^1\cdot \sum_{c_u \in B;\thinspace c \ge d} \thinspace
      q^{(a-1) + (b-2) + (c-1) + (d-1) }\thinspace
=\thinspace    
\sum_{c_u \in B;\thinspace c \ge d} \thinspace
      q^{(a-1) + (b-1) + (c-1) + (d-1) }
$
\item 
$\hspace*{0.250em}\hphantom{q^1\cdot }\sum_{c_u \in B;\thinspace c < d} \thinspace
       q^{(a-1) + (b-2) + (c-1) + (d-1) }$
\item 
$q^1\cdot \sum_{c_u \in A;\thinspace a < b} \thinspace
       q^{(a-1) + (b-2) + (c-1) + (d-1) }\thinspace
=\thinspace
\sum_{c_u \in A;\thinspace a < b} \thinspace
       q^{(a-1) + (b-1) + (c-1) + (d-1) }\thinspace .
$
\end{itemize}

(We have already proved, 
 in Section 3.1
 for 
 Garrett and Hummel's $q$-analogue,
 that 
 each of 
 these exponents of $q$ is valid as a measure of taxicab distance
 over its respective unit cubes.)

And one sees that our summation over $B_{c\ge d}$ and our summation over $A_{a<b}$
are adding the same common monomial. 
(We have already proved, 
for equation (21) in Section 3.1,
that this summation over $A_{a< b} $ and this summation over $B_{c< d}$ 
are equal as sums of their respective monomials.)
\\

We interpret the first sum of (41), using taxicab distance on the unit cubes in
$B_{c\ge d}\thinspace .$
\begin{lemma}\label{Thm:3.16}
Let $c_u$ denote an arbitrary unit cube in block $B ,$
and let the location label of $c_u$ be $(a,b,c,d) .$
For the unit cubes $c_u$ in $B$ that satisfy $c\ge d ,$
we have the equality
\begin{equation}
q^{1}\cdot \sum_{i=1}^{n} \thinspace
  [i]_q\thinspace
  q^{i-1}\thinspace
  \frac{[i+1]_q\thinspace [i]_q}{1+q}\thinspace =
  \sum_{c_u\in B;\thinspace c\ge d} 
  q^{(a-1) + (b-1) + (c-1) + (d-1)}\thinspace .
\end{equation}
The exponent of $q$ on this right side is the taxicab distance between $c_u$
and the unit cube in the region $R_3(n)$ whose location label is $(1,1,1,1).$
\end{lemma}
\begin{proof}
This is the same 
as Lemma~\ref{Thm:3.3}, 
which makes a similar statement for the unit cubes in $B_{c\ge d}$
and their taxicab distance to the start point $(1, 2, 1, 1) ,$
but 
now
with each side of equation (19)
multiplied by $q^1 .$
\end{proof}
Now we interpret the second summation of (41).
\begin{lemma}\label{Thm:3.17}
Let $c_u$ denote an arbitrary unit cube in block $B ,$
and let the location label of $c_u$ be $(a,b,c,d) .$
For the unit cubes $c_u$ in $B$ that satisfy $c <d ,$
we have the equality
\begin{equation}
q^{1}\cdot \sum_{i=1}^{n} \thinspace
  [i]_q\thinspace
  q^{i-1}\thinspace
  \frac{[i]_q\thinspace [i-1]_q}{1+q}\thinspace =
  \sum_{c_u\in B;\thinspace c < d} 
  q^{(a-1) + (b-2) + (c-1) + (d-1)}\thinspace .
\end{equation}
\end{lemma}
\begin{proof}
This is the same 
as Lemma~\ref{Thm:3.5},
which makes a similar statement for the unit cubes in $B_{c<d}$
and their taxicab distance to the start point $(1, 2, 1, 1) ,$
but 
now with each side of 
equation (20) 
multiplied by $q^1 .$
\end{proof}
Now for this subcomponent $B_{c<d}\thinspace , $ 
we apply the isometry mapping $\psi_{z,w,y-1,x}\thinspace $ from (16). 
This $\psi_{z,w,y-1,x}$ maps the point $(1,2,1,1)$
to $(1,1,1,1) .$
By the same reasoning as we used to develop Lemma~\ref{Thm:3.6}, 
we obtain a new interpretation of the second summation of (41).
\begin{lemma}\label{Thm:3.18}
Let $c_u$ denote an arbitrary unit cube in block $A ,$
and let the location label of $c_u$ be $(a,b,c,d) .$
For the unit cubes $c_u$ in $A$ that satisfy $a < b ,$
we have the equality
\begin{equation}
q^{1}\cdot \sum_{i=1}^{n} \thinspace
  [i]_q\thinspace
  q^{i-1}\thinspace
  \frac{[i]_q\thinspace [i-1]_q}{1+q}\thinspace =
  \sum_{c_u\in A;\thinspace a < b} 
  q^{(a-1) + (b-1) + (c-1) + (d-1)}\thinspace .
\end{equation}
The exponent of $q$ on this right side is the taxicab distance between $c_u$
and the unit cube in the region $R_3(n)$ whose location label is $(1,1,1,1).$
\end{lemma}
Lemmas~\ref{Thm:3.16} and~\ref{Thm:3.18} 
each give an interpretation 
of one of the two summations 
of the first parenthesized expression
on the left side of equation (41).
We conclude by interpreting their sum, and demonstrate 
it equals the first part of the right side of (41).
\begin{lemma}\label{Thm:3.19}
Let $c_u$ denote an arbitrary unit cube sized location in the region $R_3(n) ,$
and let the location label of $c_u$ be $(x,y,z,w) = (a,b,c,d) .$
We have the equalities
\begin{equation*}
\begin{split}
\sum_{c_u \in B;\thinspace c \ge d}
 q^{(a-1) + (b-1) + (c-1) + (d-1)}\thinspace
+ 
\sum_{c_u \in A;\thinspace a < b}
 q^{(a-1) + (b-1) + (c-1) + (d-1)}
&=\\
\sum_{\substack{
      c_u \in R_3(n)\\
      1 \le a <   b \le n+1\\ 
      1 \le d \le c \le n
                }
     }
 q^{(a-1) + (b-1) + (c-1) + (d-1)}
&=\thinspace
q^1\cdot {\genfrac{[}{]}{0pt}{0}{n+1}{2}}_q^2 \thinspace .
\end{split}
\end{equation*}
\end{lemma}
\begin{proof}
This is the same 
as Lemma~\ref{Thm:3.7},
which makes a similar statement for the unit cubes in $S_{3,3}(n) $
and their taxicab distance to the start point $(1,2,1,1) ,$
but
now
with each side of its equations 
multiplied by $q^1 .$
\end{proof}
Here we begin step 2 of interpretation of equation (41):
for the second parenthesized expression on the left side,
incorporating its outside factor $q^2 .$
\\

Select the block $B.$
Let $c_u$ denote an arbitrary unit cube in $B,$ and let its location label
be $(x,y,z,w) = (a,b,c,d) .$ 
Partition $B$ into the parts $B_{c\ge d}$ and $B_{c< d}\thinspace .$

Introduce a new isometry mapping, 
$\psi_{y, x+1}\thinspace ,$ defined by the rule
\begin{equation*}
\psi_{y, x+1}\thinspace = \tau_{y+1} \circ \phi_{y,x} 
\end{equation*}
where
$\phi_{y,x}: (x,y,z,w) \mapsto (y,x,z,w)$ is the reflection 
that permutes the axes $x$ and $y$; and
$\tau_{y+1}: (x,y,z,w) \mapsto (x, y+1, z, w)$ is 
translation by $+1$ in $y .$ 

One checks that this mapping $\psi_{y, x+1}$ 
acts on each set of unit cubes that we 
had used in our previous construction, 
for interpreting the first parenthesized expression on the left side of (41),
by sending it to the following image set:
\begin{itemize}
\item
  $\psi_{y, x+1}$ sends the block $B$ to block $C ,$ 
   for each block in its defined location.
\item  
  $\psi_{y, x+1}$ sends the partition part $B_{c \ge d}\thinspace ,$
  to the subcomponent of the block $C$ which comprises the unit cubes $c_u$
  whose location labels $(a,b,c,d)$ satisfy $c \ge d .$
  
  And $\psi_{y, x+1}$ sends the partition part $B_{c < d}\thinspace ,$
  to the subcomponent of $C$ comprising the unit cubes which satisfy $c < d .$
  
   These image sets form a partition of $C .$
   We call them respectively, $C_{c \ge d}$ and $C_{c < d} .$
\item  
  $\psi_{y, x+1}$ sends the part $A_{a < b}\thinspace ,$
  to the subcomponent of block $A$ 
  which comprises the unit cubes whose location labels satisfy 
  $2 \le b \le a \le n .$
  We call this image set $A_{2 \le b \le  a} .$
\item    
  $\psi_{y, x+1}$ sends the union of the disjoint sets of unit cubes,
  $B_{c\ge d}\thinspace \bigcup A_{a < b} = S_{3,3}(n),$ 
  to the union of the disjoint sets of unit cubes,
  $C_{c \ge d}$ and $A_{2 \le b \le a} \thinspace .$
  
  These we assemble to form a new solid in the region $R_3(n) ,$
  which comprises the unit cubes of $R_3(n)$ that satisfy both 
  $2 \le b \le a \le n+1$ and $1 \le d \le c \le n .$
  In other words they form the image of the set $S_{3,3}(n) ,$
  which by definition satisfies $1 \le a < b \le n+1 ,$
  under the isometry mapping $\psi_{y,x+1}\thinspace .$
  
  We will denote this new solid by 
  $\thinspace C_{c\ge d}\thinspace \bigcup A_{2 \le b \le a} 
  \thinspace
  =\thinspace 
  S_{3,3}^{\prime}(n)\thinspace .$
\end{itemize}  

Now we apply these images under the isometry mapping $\psi_{y,x+1}\thinspace ,$
to interpret the second parenthesized expression on the left side of (41).
\\

For this $q$-analogue equation, we will select for our $q$-counts of 
taxicab distances:
\begin{itemize}
\item 
$
q^1 \cdot q^1 \cdot \sum_{c_u \in C;\thinspace c \ge d} \thinspace
      q^{(a-2) + (b-2) + (c-1) + (d-1) }\thinspace
=\thinspace
\sum_{c_u \in C;\thinspace c \ge d} \thinspace
      q^{(a-1) + (b-1) + (c-1) + (d-1) }
$
\item 
$
\hspace*{0.37500em}\hphantom{q^1 \cdot }q^1 \cdot \sum_{c_u \in C;\thinspace c < d} \thinspace
       q^{(a-2) + (b-2) + (c-1) + (d-1) }\thinspace
=\thinspace
\sum_{c_u \in C;\thinspace c < d} \thinspace
       q^{(a-1) + (b-2) + (c-1) + (d-1) }
$
\item 
$
q^1\cdot  q^1\cdot \sum_{c_u \in A;\thinspace 2 \le b \le a } \thinspace
       q^{(a-2) + (b-2) + (c-1) + (d-1) }\thinspace
=\thinspace
\sum_{c_u \in A;\thinspace 2 \le b \le a } \thinspace
       q^{(a-1) + (b-1) + (c-1) + (d-1) }\thinspace .
$
\end{itemize}

One checks that each of these exponents of $q$ is valid as a
measure of taxicab distance over its respective unit cubes.

For 
any unit cube $c_u$ in $C_{c \ge d}\thinspace ,$ 
it is valid to measure a taxicab distance by using
$2 + (a-2) + (b-2)  + (c-1) + (d-1)\thinspace $:
its 
coordinate values $a,b$ are each at least 2;
and its coordinates $c,d$ are each at least 1.

For 
any unit cube $c_u$ in $C_{c < d}\thinspace ,$ 
it is valid to measure a taxicab distance by using
$1 + (a-2) + (b-2)  + (c-1) + (d-1)\thinspace $:
its 
coordinate values $a,b,d$ are each at least 2;
and its 
coordinate $c$ is at least 1.

For 
any unit cube $c_u$ in $A_{2 \le b \le a }\thinspace ,$
it is valid to measure a taxicab distance by using
$2 + (a-2) + (b-2)  + (c-1) + (d-1)\thinspace $:
its 
coordinate values $a,b$ are each at least 2;
and its coordinates $c,d$ are each at least 1.
\\

And one sees that our summation over $C_{c\ge d}$ and our summation over $A_{2\le b \le a}$
are adding the same common monomial. We will prove ahead, 
using the composite isometry map 
$\psi_{y, x+1} \circ \psi_{z, w, y-1, x}\thinspace, $ 
that the summation over $A_{2 \le b\le a} $ and the
summation over $C_{c< d}$ are equal as sums of their respective monomials.
\\

We interpret the first sum,
of the second parenthesized expression on the left side of (41),
using taxicab distance on the unit cubes in 
$\thinspace \psi_{y,x+1}(B_{c\ge d})=  C_{c\ge d}\thinspace .$

\begin{lemma}\label{Thm:3.20}
Let $c_u$ denote an arbitrary unit cube in block $C ,$
and let the location label of $c_u$ be $(a,b,c,d) .$
For the unit cubes $c_u$ in $C$ that satisfy $c\ge d ,$
we have the equality
\begin{equation}
q^{2}\cdot \sum_{i=1}^{n} \thinspace
  [i]_q\thinspace
  q^{i-1}
  \frac{[i+1]_q\thinspace [i]_q}{1+q}\thinspace =
  \sum_{c_u\in C;\thinspace c\ge d} 
  q^{(a-1) + (b-1) + (c-1) + (d-1)}\thinspace .
\end{equation}
The exponent of $q$ on this right side is the taxicab distance between $c_u$
and the unit cube in the region $R_3(n)$ whose location label is $(1,1,1,1).$
\end{lemma}
\begin{proof}
Begin with Lemma~\ref{Thm:3.16}, 
which makes a similar statement for 
the unit cubes in $B_{c\ge d}$
and their taxicab distance to the start point $(1, 1, 1, 1) .$
Now recall that applying the isometry mapping $\psi_{y,x+1}\thinspace $
to the unit cubes in $B_{c\ge d} \thinspace $ 
will send them onto the unit cubes in $C_{c\ge d} \thinspace .$
And one checks that $\psi_{y, x+1}$ sends the point $(1,1,1,1)$ to $(1,2,1,1) .$
So we obtain the equalities
\begin{equation*}
\begin{split}
  q^{1}\cdot \sum_{c_u\in B;\thinspace c \ge d} 
  q^{(a-1) + (b-1) + (c-1) + (d-1)}\thinspace 
=\thinspace
  q^{1}\cdot 
  &\sum_{c_u\in C;\thinspace c \ge d} 
  q^{(a-1) + (b-2) + (c-1) + (d-1)}\\
=\thinspace
  &\sum_{c_u\in C;\thinspace c \ge d} 
  q^{(a-1) + (b-1) + (c-1) + (d-1)}\thinspace .  
\end{split}
\end{equation*}
Now apply the equation statement (42) from Lemma~\ref{Thm:3.16}, 
to complete our proof. 
\end{proof}  

Now we interpret the second summation,
of the second parenthesized expression on the left side of (41),
using the remaining part of our partition, 
$\thinspace \psi_{y,x+1}(B_{c < d}) =  C_{c < d}\thinspace .$

\begin{lemma}\label{Thm:3.21}
Let $c_u$ denote an arbitrary unit cube in block $C ,$
and let the location label of $c_u$ be $(a,b,c,d) .$
For the unit cubes $c_u$ in $C$ that satisfy $c <d ,$
we have the equality
\begin{equation}
q^{2}\cdot \sum_{i=1}^{n} \thinspace
  [i]_q\thinspace
  q^{i-1}
  \frac{[i]_q\thinspace [i-1]_q}{1+q}\thinspace =
  \sum_{c_u\in C;\thinspace c < d} 
  q^{(a-1) + (b-2) + (c-1) + (d-1)}\thinspace .
\end{equation}
\end{lemma}
\begin{proof}
Begin with Lemma~\ref{Thm:3.17}, 
which makes a similar statement for the unit cubes in $B_{c < d}$
and their taxicab distance to the start point $(1,2,1,1) .$
Now recall that applying the isometry mapping $\psi_{y, x+1}$
to the unit cubes in $B_{c < d} $
will send them onto the unit cubes in $C_{c < d}\thinspace .$
And one checks that $\psi_{y, x+1}$ sends the point $(1,2,1,1)$ to $(2,2,1,1).$
So we obtain the equalities
\begin{equation*}
\begin{split}
  q^{1}\cdot \sum_{c_u\in B;\thinspace c < d} 
  q^{(a-1) + (b-2) + (c-1) + (d-1)}\thinspace 
=\thinspace
  q^{1}\cdot 
  &\sum_{c_u\in C;\thinspace c < d} 
  q^{(a-2) + (b-2) + (c-1) + (d-1)}\\
=\thinspace
  &\sum_{c_u\in C;\thinspace c < d} 
  q^{(a-1) + (b-2) + (c-1) + (d-1)}\thinspace .  
\end{split}
\end{equation*}
Now apply the equation statement (43) from Lemma~\ref{Thm:3.17}, 
to complete our proof.
\end{proof}

Next we have 
new interpretation of the second summation,
of the second parenthesized expression on left side of (41),
by using the image of the (composite) isometry mapping, 
$\thinspace 
  \psi_{y, x+1}(\psi_{z,w, y-1, x} (B_{c<d}) )\thinspace =\thinspace 
  \psi_{y, x+1}(A_{a<b})\thinspace =\thinspace 
  A_{2 \le b \le a}\thinspace .$
\begin{lemma}\label{Thm:3.22}
Let $c_u$ denote an arbitrary unit cube in block $A ,$
and let the location label of $c_u$ be $(a,b,c,d) .$
For the unit cubes $c_u$ in $A$ that satisfy $2\le b \le  a ,$
we have equality
\begin{equation}
q^{2}\cdot \sum_{i=1}^{n} \thinspace
  [i]_q\thinspace
  q^{i-1}
  \frac{[i]_q\thinspace [i-1]_q}{1+q}\thinspace =
  \sum_{c_u\in A;\thinspace 2 \le b \le a} 
  q^{(a-1) + (b-1) + (c-1) + (d-1)}\thinspace .
\end{equation}
The exponent of $q$ on this right side is the taxicab distance between $c_u$
and the unit cube in the region $R_3(n)$ whose location label is $(1,1,1,1).$
\end{lemma}
\begin{proof}
Begin with Lemma~\ref{Thm:3.18}, 
which makes a similar statement for the unit cubes 
in $A_{a < b}$
and their taxicab distance to the start point $(1,1,1,1) .$
Now recall that applying the isometry mapping $\psi_{y, x+1}$
to the unit cubes in $A_{a < b} $
will send them onto the unit cubes in $A_{2 \le b \le a}\thinspace .$
And 
$\psi_{y, x+1}$ sends the point $(1,1,1,1)$ to $(1,2,1,1).$
So we obtain the equalities
\begin{equation*}
\begin{split}
  q^{1}\cdot \sum_{c_u\in A;\thinspace a < b} 
  q^{(a-1) + (b-1) + (c-1) + (d-1)}\thinspace 
=\thinspace
  q^{1}\cdot 
  &\sum_{c_u\in A;\thinspace 2 \le b \le a } 
  q^{(a-1) + (b-2) + (c-1) + (d-1)}\\
=\thinspace
  &\sum_{c_u\in A;\thinspace 2 \le b \le a } 
  q^{(a-1) + (b-1) + (c-1) + (d-1)}\thinspace .  
\end{split}
\end{equation*}
Now apply the equation statement (44) from Lemma~\ref{Thm:3.18}, 
to complete our proof.
\end{proof}

Lemmas~\ref{Thm:3.20} and~\ref{Thm:3.22} 
each give an interpretation of one of the two summations,
of the second parenthesized expression 
on the left side of equation (41).
We conclude by interpreting their sum, and demonstrate 
it equals the second part of right side of (41).

\begin{lemma}\label{Thm:3.23}
Let $c_u$ denote an arbitrary unit cube sized location in the region $R_3(n) ,$
and let the location label of $c_u$ be $(x,y,z,w) = (a,b,c,d) .$
We have the equalities
\begin{equation*}
\begin{split}
\sum_{c_u \in C;\thinspace c \ge d}
 q^{(a-1) + (b-1) + (c-1) + (d-1)}\thinspace
+ 
\sum_{c_u \in A;\thinspace 2 \le b \le a}
 q^{(a-1) + (b-1) + (c-1) + (d-1)}
&=\\
\sum_{\substack{
      c_u \in R_3(n)\\
      2 \le b \le   a \le n+1\\ 
      1 \le d \le c \le n
                }
     }
 q^{(a-1) + (b-1) + (c-1) + (d-1)}
&=\thinspace
q^2\cdot {\genfrac{[}{]}{0pt}{0}{n+1}{2}}_q^2 \thinspace .
\end{split}
\end{equation*}
\end{lemma}
\begin{proof}
The first equation follows from our discussion 
of the isometry $\psi_{y,x+1}\thinspace ,$
that the union of the disjoint sets of unit cubes 
$C_{c\ge d}$ and $A_{2 \le b \le a}$ is the set
$S_{3,3}^{\prime}(n)$ of unit cubes in $R_3(n)$ which satisfy
both $2 \le b \le a \le n+1$ and $1 \le d \le c \le n .$

To establish the second equation,
recall 
from our proof of Lemma~\ref{Thm:3.20}
and our proof of Lemma~\ref{Thm:3.22} 
that applying $\psi_{y,x+1}$
to the unit cubes in $B_{c \ge d}$ 
and in $A_{a<b}$
will increase each of their taxicab distances to $(1, 1, 1, 1)$ by $+1 .$
Now apply the equation 
from Lemma~\ref{Thm:3.19}, 
which makes a similar statement for the unit cubes in $S_{3,3}(n)$
and their taxicab distance to the start point $(1,1,1,1),$
to complete our proof.
\end{proof}

Here is step 3 of our interpretation of the $q$-analogue equation (41),
for the total: 
of the first parenthesized expression on the left side
(incorporating its outside factor $q^1$),
plus  
second parenthesized expression on the left side
(incorporating outside factor $q^2$).
We will interpret by using 
the union of these four disjoint sets of unit cubes,
to form a reassembled solid in $R_3(n).$

\begin{lemma}\label{Thm:3.24}
Let $c_u$ denote an arbitrary unit cube sized location in the region $R_3(n) ,$
and let the location label of $c_u$ be $(x,y,z,w) = (a,b,c,d) .$
We have the equalities
\begin{equation*}
\begin{split}
\sum_{\substack{
      c_u \in R_3(n)\\
      1 \le a <   b \le n+1\\ 
      1 \le d \le c \le n
                }
     }
 q^{(a-1) + (b-1) + (c-1) + (d-1)}\thinspace
&+ \thinspace
\sum_{\substack{
      c_u \in R_3(n)\\
      2 \le b \le   a \le n+1\\ 
      1 \le d \le c \le n
                }
     }
 q^{(a-1) + (b-1) + (c-1) + (d-1)}
\\
&= \sum_{\substack{
      c_u \in R_3(n)\\
      1 \le a \le n+1;\thinspace
      2 \le b \le n+1\\ 
      1 \le d \le c \le n
                }
     }
 q^{(a-1) + (b-1) + (c-1) + (d-1)}\thinspace
=\thinspace
(q + q^2)\thinspace {\genfrac{[}{]}{0pt}{0}{n+1}{2}}_q^2 \thinspace .
\end{split}
\end{equation*}
\end{lemma}
\begin{proof}
The first equation follows from two observations.

These two sets of unit cubes over which we are summing, 
$S_{3,3}(n)$ and $S_{3,3}^{\prime}(n) ,$ are disjoint: 
Each unit cube $c_u$ in the first set
has a location label $(a,b,c,d)$ which satisfies $a < b ;$
and each unit cube in the second set must satisfy $b \le a .$

The union of these two sets 
comprises all the unit cubes in $R_3(n)$ 
whose location labels obey the constraints that:
the coordinate values $(c,d)$ satisfy $1 \le d \le c \le n ;$
the coordinate value $b$ satisfies $2 \le b \le n+1 ;$ and 
the coordinate value $a$ satisfies $1 \le a < b $ or $b \le a \le n+1 .$
So the coordinate value $a$ has no constraints 
within its domain $\{1, 2, \ldots, n+1\}$ in $R_3(n) .$%

And now the second equation follows by applying Lemma~\ref{Thm:3.19} 
and Lemma~\ref{Thm:3.23}; 
the sum of the right sides of their equations is our total. 
\end{proof}
In Lemma~\ref{Thm:3.24},
our final assembly $S_{3,3}(n) \bigcup S_{3,3}^{\prime}(n)$
is the set of unit cubes $c_u$ in $R_3(n)$
whose location label coordinates $(a,b,c,d)$ satisfy that
$1 \le a \le n+1 , $
$2 \le b \le n+1 ,$ and
$1 \le d \le c \le n .$
The shape of this assembled solid
is an $(n+1) \times n$ rectangle in the $(x,y)$-plane,
for each value of $(z,w)$ in $R_3(n).$
And it is an (oblique) step triangle in the $(z,w)$-plane,
for each value of $(x,y)$ in $R_3(n) .$

For nomenclature, we will call the shape of
this 4-dimensional assembled solid an
\emph{(oblique) step triangular -- rectangular duoprism}.
\\

This completes our interpretation of the $q$-analogue equation from 
Forster, Garrett, Jacobsen, and Wood, as stated in the form (41),
and omitting the outside factor 2.

\section{Further work}

For higher degree power sums, 
A. F. Beardon~\cite{example8_beardon}  
presented general identities using 
``ideas [that] date back to Johann Faulhaber (1631)".
Let ${\sigma_p(n) = 1^p + 2^p + \ldots + n^p .}$
Beardon's identities state, in a form we have mildly rewritten, 
that for each $p \ge 3$:

\begin{equation*}
\begin{split}
\textrm{if }p\textrm{ is odd:}\hspace*{2.320em}  
n^{\frac{p+1}{2}} (n+1)^{\frac{p+1}{2}}
&=
\sum_{\substack{1 \le r \le {\frac{p+1}{2}}\\
                r\textrm{ is odd} } } 
     2\thinspace \binom{\frac{p+1}{2}}{r}\thinspace \sigma_{p+1-r}(n)\\
\textrm{if }p\textrm{ is even:}\medspace
(n + \tfrac{1}{2})\thinspace n^{\frac{p}{2}} (n+1)^{\frac{p}{2}}
&=
\sum_{\substack{1 \le r \le {\frac{p}{2}}\\
                r\textrm{ is odd} } } 
     2\thinspace \binom{\frac{p}{2}}{r}\thinspace \sigma_{p+1-r}(n)
     +
\sum_{\substack{0 \le r \le {\frac{p}{2}}\\
                r\textrm{ is even} } } 
     \binom{\frac{p}{2}}{r}\thinspace \sigma_{p-r}(n) \thinspace .     
\end{split}
\end{equation*}
For example, 
the identity for the sum of $4^{\textrm{th}}$ powers,
and for the sum of $5^{\textrm{th}}$ powers:
\begin{align}
(n + \tfrac{1}{2})\thinspace n^2 (n+1)^2 
&
= 
5 \sigma_4(n) + \sigma_2(n)  \\
n^3 (n+1)^3 
&= 
6 \sigma_5(n) + 2 \sigma_3(n)\thinspace . 
\end{align}
We have established 
geometric assemblies
which generalize our 
four-block construction
for $p=3$
presented in 
this article, to interpret these identities 
(48) and (49) 
and for higher powers.
Those results will be the subject of another article, 
which is in preparation.

\end{document}